\documentstyle[12pt]{article}
\begin{document}

\newtheorem{ct}{\indent Connection}
\newtheorem{th}{\indent  Theorem}
\newtheorem{df}{\indent  Definition}
\newtheorem{prop}{\indent Proposition}
\newtheorem{rk}{\indent Remark}
\newtheorem{lm}{\indent Lemma}
\newtheorem{prb}{\indent Problem}
\newtheorem{con}{\indent Conjecture}
\newtheorem{cy}{\indent Corollary}

\title{ Hyperbolic groups and free constructions}
\author{O. Kharlampovich, A. Myasnikov}
\maketitle
Hyperbolic groups have been the subject of intensive investigation since the 
work 
of 
Gromov \cite{G}.  
Let $G=<J,\cal R>$ be a finitely presented group with a set of
generators $J$ and a set of relators $\cal R.$ A word $W$ in the
alphabet $J^{\pm 1}$ is equal to $1$ in $G$ if and only if there is
an equality 
\begin{equation}
\label{hyp}
W={\prod}_{i=1}^n S_i^{-1}R_i^{\pm 1}S_i\end{equation}
 in the free group $F=F(J),$
where $S_i\in F$ and $R_i\in \cal R.$ The group $G$ is {\em hyperbolic}
if there exists a linear function bounding the minimal number
of factors $n=n(W)$ in (\ref{hyp}) 
depending on the length $\mid\mid W\mid\mid $ of the
word $W.$ This definition does not depend on the choice of the
presentation of $G.$

In his book (\cite{G}, p.3.3) 
M.Gromov claimed that if $G_1$ and $G_2$ are torsion-free 
hyperbolic and $U$ and $V$ are maximal cyclic subgroups in $G_1$ and 
 $G_2$ respectively, then the amalgamated free product 
${G}_1\ast_{U=V}{G}_2$ is also hyperbolic.

On the other hand the  group $<x>\ast_{x^n=y^m}<y>$
for $\mid n\mid ,\mid m\mid >1$ is not hyperbolic, because it contains
a free abelian subgroup with generators $xy$ and $x^n.$

The groups of Baumslag-Solitar BS$(m,n)= <x,t\mid t^{-1}x^mt=x^n>$
provide examples of HNN-extensions of hyperbolic group 
(with cyclic associated subgroups) that are not hyperbolic.

In \cite{BGSS} it has been shown that an amalgamated product of two hyperbolic
groups with a cyclic subgroup amalgamated is automatic and that an amalgamated 
product of two finitely generated free groups with a finitely generated 
subgroup amalgamated is acynchronously automatic. It has been also proved
in \cite{BGSS} that if $G$ is an HNN-extension of a finitely generated free 
group
with finitely many stable letters and if the associated subgroups are all
finitely generated, then $G$ is asynchronously automatic. Under the
additional assumption of ``speed-matching'', such HNN-extensions are
shown to be automatic in \cite{Sh}. (Notice that the 
class of automatic groups is contained in the class of asynchronously
automatic groups and contains the class of hyperbolic groups.) 
In \cite{BGSS1}, \cite{GS} it was proved that if an amalgamated free product
${G}_1\ast_{U}{G}_2$
is automatic then both groups $G_1$ and $G_2$ 
are automatic provided the amalgamated
subgroup $U$ is finite.

In \cite{BF} Bestvina and Feighn proved the combination theorem for negatively 
curved spaces and as a corollary obtained the result which we also formulate
in this paper as  Corollary 2.

We define
a subgroup ${U}$ of a group $G$ to be {\em conjugate separated} if
the set $\{u\in U \mid u^x\in U\}$ is finite for all $x \in G \setminus U$.

Let us introduce the construction of a separated HNN-extension of
a group $G$. 
\begin{df}
Suppose that $U$ and $V$ are  subgroups of $G$,
$\psi: U \rightarrow V$ is an isomorphism,
either $U$ or $V$ 
is conjugate separated, and the set $U \cap g^{-1}Vg$ is finite 
for all $g \in G$. 
Then the HNN-extension $\left< G, t \mid t^{-1}ut=u^\psi, \ u \in U
\right>$ is called separated. 
\end{df}

Let $a,b\in G$; by $\left |a-b\right |_G$ we denote the distance between  
the points $a$ and $b$ in the Cayley graph of $G$ (see the definition below).
If it is clear from the context in which group the distance is taken,
we will just write $\left |a-b\right |.$

A finitely generated subgroup $U$ of a hyperbolic group $G$ is said
to be 
{\em quasiisometrically embedded} 
if there is a constant $\lambda =\lambda (U)$ such that
$\left |1-a\right |_U\leq \lambda \left |1-a\right |_G$ for any 
element $a\in U.$ 
Every quasiisometrically embedded subgroup of a hyperbolic group
is itself hyperbolic. 

It is not hard to see that a  
subgroup of a hyperbolic group is 
quasiisometrically embedded if and only if it is quasiconvex
(the definition of quasiconvexity can be found in Section 6).
\begin{th} If ${\cal G}$ is a hyperbolic group
and ${\cal H}=<{\cal G},t\mid  {\cal U}^t={\cal V}>$
is a separated
HNN-extension such that the subgroups ${\cal U}$ and ${\cal V}$ 
are quasiisometrically embedded in ${\cal G},$ then
${\cal H}$ is hyperbolic.\end{th}

\begin{th} 
Let ${\cal G}_1$, ${\cal G}_2$ be hyperbolic groups, 
${\cal U}\leq {\cal G}_1$, ${\cal V}\leq 
{\cal G}_2$ quasiisometrically embedded, and  ${\cal U}$ conjugate separated
in ${\cal G}_1.$
Then the group ${\cal G}_1\ast_{{\cal U}={\cal V}}{\cal G}_2$
is hyperbolic.
\end{th}

  As corollaries we have the following results.

\begin{cy} If ${\cal G}$ is a hyperbolic group,
$\cal A$ and $\cal B$ isomorphic virtually cyclic subgroups 
, then the HNN-extension
${\cal H}=<{\cal G},t\mid  {\cal A}^t={\cal B}>$
is hyperbolic if and only if it is separated.
\end{cy}

\begin{cy} (\cite{BF}) Let ${\cal G}_1$, ${\cal G}_2$ be hyperbolic groups, 
${\cal A}\leq {\cal G}_1$, ${\cal B}\leq 
{\cal G}_2$, virtually cyclic.  
Then the group ${\cal G}_1\ast_{{\cal A}={\cal B}}{\cal G}_2$
is hyperbolic if and only if either ${\cal A}$ is conjugate separated
in ${\cal G}_1$ or ${\cal B}$ is conjugate separated
in ${\cal G}_2.$
\end{cy}

The assertion of our Corollary 1  contradicts the assertion of the
last corollary (HNNs over virtually cyclics) in \cite{BF} (there is
an omission in this corollary in \cite{BF}). The group
given by the presentation $K=<a,b,t\mid t^{-1}a^2t=b^2>$ is a counter-example.
The group is obviously 
not hyperbolic, but it satisfies condition (2) in the last corollary in
\cite{BF}. 

After our paper was ready, A.Yu. Olshanskii informed us that he and
his student K. Mikhailovskii \cite{OM} independently obtained the results
formulated in Corollaries 1 and 2. Corollary 2 is also proved in
\cite{Rita}. 

\begin{cy}
HNN-extensions (amalgamated products) of hyperbolic groups
with finite associated (amalgamated) subgroups are hyperbolic.\end{cy}

\begin{cy} Separated HNN-extensions of a free group with finitely generated
associated subgroups are hyperbolic.
\end{cy}
The condition in Theorems 1 and 2 can be weakened (see Section 5).
The proof of our results  
uses a geometric interpretation by 
 Van-Kampen diagrams of the
deducibility of relations
in a group from the defining relations, as well as the hyperbolicity
of the Cayley graph $\Gamma ({\cal G})$ of a hyperbolic group ${\cal G}.$

In Section 6 we prove some results on quasiconvexity (Theorems 4, 5, 6 and 7).

In Section 7 we apply our results to exponential groups.

Let $A$ be an arbitrary associative ring with identity and $G$ a group.
Fix an action of the ring $A$ on $G$, i.e. a map $G \times A \rightarrow
G$. The result of the action of $\alpha \in A$ on $g \in G$ is written
as $g^\alpha$.
Groups with $A$-actions satisfying axioms 1)--4) in Section 7 are called  
{\em $A$--groups}.   
In particular, an arbitrary group $G$ is a ${\bf Z}$-group.
In the case where $A$ is the field 
of rationals ${\bf Q},$  $\bf Q$-completions of groups (see definition
in
Section 7) were studied by G. Baumslag in
\cite{bau}, \cite{baum1}. $A$-completions for arbitrary rings $A$
were investigated in \cite{MyasExpo2}. 

In \cite{baum1} it was proved that for a free group $F$ 
the word problem in the  $\bf Q$-completion $F^{\bf Q}$
is solvable. The proof was based on the residual finiteness
of some specific subgroups of $F^{\bf Q}$, so the problem was posed 
of finding a ``direct'' proof of the solvability of the 
word problem using normal forms of 
elements in $F^{\bf Q}$. Moreover, in the same article G.Baumslag 
mentioned that the conjugacy problem in $F^{\bf Q}$ is still  open.

   In Section 7 we describe the ${\bf Q}$-completion $G^{\bf Q}$ 
of a torsion-free
hyperbolic group $G$ as the union of an effective chain of hyperbolic subgroups.
This allows one to apply techniques of hyperbolic group
theory  to solve various algorithmic
problems in $G^{\bf Q}$, in particular, to solve the conjugacy
problem (Theorem 10) and to construct effectively 
some natural  normal forms for its elements, induced
by the normal forms of amalgamated 
free products (Theorem 9). 
A free group $F$ is hyperbolic, so one can answer two of
Baumslag's questions \cite{baum1} mentioned above.

\section{Quasigeodesic Polygons in Hyperbolic Groups}

Let us recall some notions from the theory of hyperbolic spaces.

Let $X$ be a metric space, $|x-y|$ the distance between points $x,y\in X.$
If one fixes a point $o\in X$ then Gromov's product $(x\cdot y)_o$ 
is by definition
$$(x\cdot y)_o=1/2(|x-o|+|y-o|-|x-y|).$$
The space $X$ is called $\delta$-hyperbolic (for some fixed constant $\delta
\geq 0$) if for all points $x,y,z,o\in X$
$$(x\cdot y)_o\geq \min ((x\cdot z)_o,(y\cdot z)_o)-\delta .$$
And $X$ is called hyperbolic if it is $\delta$-hyperbolic
for some $\delta \geq 0$

A geodesic segment between points $x,y\in X$ is an isometric map 
$[0,|x-y|]\longrightarrow X$ 
sending $0$ to $x$ and $|x-y|$ to $y.$ Its image will also be called
a geodesic segment;
$[x,y]$ is the notation for some fixed geodesic segment between $x$ and
$y.$ A metric space is called geodesic if every pair of points
can be connected by a geodesic segment. We shall call a geodesic 
$\delta$-hyperbolic space simply a $\delta$-space.

An example of a geodesic space is the realization of the Cayley graph $\Gamma 
({G})=
\Gamma ({G},J)$ of a group ${G}$ with a fixed generating system $J.$
Recall that the vertices of $\Gamma ({G})$ are elements of ${G},$ and
the elements $g,\ h=ga,\ a\in J, $ are connected by an edge $e=(g,a)$ having 
label
$\phi (e)=a\in J.$ The label of a path is the product of the labels of
the edges of this path. Let us endow each edge $e$ with the metric of the unit 
segment $[0,1].$ By definition now the distance $\left | x-y\right | $
between points $x$ and $y$ is the length of a shortest path in 
$\Gamma ({G})$ between $x$ and $y.$

A finitely presented 
group ${G}$ is hyperbolic if and only if $\Gamma ({G})$ is a hyperbolic
space \cite{G}. If $g,h\in G$ then by definition 
$(g\cdot h)=(g\cdot h)_1$ is the
Gromov product in $\Gamma ({G}),$ where $1$ is the identity of ${G}.$
Let $p$ be a path and $x$ a point. 
Let $\left |x,p\right |$ denote the distance between the point $x$ and the 
path $p$ and $\left |x,[y,z]\right |$ denote the distance between the point $x$ 
and the 
geodesic segment
$[y,z].$
\begin{lm} (\cite{Ol}, lemma 1.5) 
\label{1.5}
For each geodesic triangle $[x_1,x_2,x_3]$ in a $\delta$-space,
there are points $y_i\in [x_{i-1},x_{i+1}]$ (indices are considered
modulo 3) such that
$$\left | x_i-y_{i-1}\right | =\left | x_i-y_{i+1}\right | =(x_{i-1}\cdot 
x_{i+1})_{x_i},$$
$$\left | y_i-y_{i-1}\right |\leq 4\delta \ \ and 
\ \left |u,[x_i,y_{i\pm 1}]\right |\leq 4\delta $$
for any point $u\in [x_i,y_{i\pm 1}].$
\end{lm}

It is easy to verify that for a hyperbolic group $G$ we have 
in $\Gamma ({G})$
$$ (x\cdot z)_y=(\phi ([y,x])\cdot \phi ([y,z])).$$ We can rewrite
the equalities from the Lemma in the following form.
$$\left | x_i-y_{i-1}\right | =\left | x_i-y_{i+1}\right | =
(\phi [x_{i},x_{i-1}]\cdot \phi [x_i,x_{i+1}]).$$

A path $p$ with the natural parametrisation by length in $\Gamma ({G})$
is called {\em $(\lambda ,\mu )$-quasigeodesic} for some  $\lambda\geq 0$
and $\mu\geq 0$, if for any points $p(s)$ and $p(t)$
$$\lambda \left | s-t\right |-\mu\leq \left | p(s)-p(t)\right | .$$

Denote by $q_{-}\ (resp. q_{+})$ the initial (resp. terminal) vertex of a path 
$q.$
A word in the generators of $G$ is called {\em geodesic (resp. quasigeodesic)} 
if
the corresponding path is geodesic (resp. quasigeodesic) in the Cayley graph of 
$G.$

If a subgroup $U$ is quasiisometrically embedded in
a hyperbolic group $G$ then every
geodesic in $U$ is a $(\lambda ,0)$-quasigeodesic in $G.$
\begin{lm} (\cite{Ol}, lemma 1.9)
\label{1.9}
There exists a constant $H_1=H_1(\delta ,\lambda ,\mu)$ such that for any
$(\lambda , \mu)$-quasigeodesic path $p$ in a $\delta$-space and any
geodesic path $q$ with the conditions $p_-=q_-$ and $p_+=q_+,$ the 
inequalities $\left |u,p\right |<H_1$ and $\left |v,q\right |<H_1$ 
hold for any points
$u\in q$ and $v\in p.$
\end{lm}

As in \cite{Ol}, call two paths $p$ and $q$ {\em $C$-bound}  if
$(|p_- -q_-|, |p_+ -q_+|)\leq C.$
\begin{lm}
\label{1.7}
(\cite{Ol}, lemma 1.7)Let $[x_1,\ldots ,x_4]$ be a geodesic quadrangle in a 
$\delta$-space
and
$$\left | x_1-x_2\right | >4\max(\left | x_1-x_4\right | ,\left | x_2-x_3\right 
| 
).$$
Then, there exist $8\delta$-bound subsegments $p$ and $q$ of the segments
$[x_1,x_2]$ and $[x_3,x_4]$ such that
$$min(\left | p\right | ,\left | q\right | )\geq (7/20)\left | x_1-x_2\right | 
-8\delta .$$
\end{lm}
Now we fix some notation and introduce
constants $\delta ,\ \lambda ,\ c,\ H_1,$ that will be used below
without reference. Let $\cal G$ be the hyperbolic group from 
Theorem 1,
$J$ a fixed set of generators of $\cal G.$ Suppose $G$ is
$\delta$-hyperbolic.
Let $
\cal U$ and $\cal V$ be the subgroups of $\cal G$ from Theorem 1,
quasiisometrically embedded with the constant $\lambda .$
Let $\cal U$ be conjugate separated.

We introduce two length functions on the group ${\cal G}.$
If $g\in {\cal G}$ then $\left | g\right | =\left | g-1\right |$ 
in $\Gamma ({\cal G},J).$
In other words, $\left | g\right |$ is the length of a shortest word in the
alphabet $J$
representing $g$ in ${\cal G}.$
We also consider words in the alphabet $J$ to be elements of ${\cal G}.$
The notation $W\equiv V$ means the equality of words, and $W=V$ means
the equality of elements.  Suppose $V$ is a word in the alphabet $J;$
we say that $W$ is a {\em geodesic
word} such that $W=V$ if $W$ is a shortest word representing
in $\cal G$ the same element as $V.$

By $\left |\left | W\right |\right |=\left |\left | W\right |\right
|_{\cal G}$ we denote the length of the word $W.$
Clearly the length $\left |\left | W\right |\right |$ can be greater then 
 $\left | W\right | .$

We now fix some presentation of $\cal H.$
Let $A=\{a_1,\ldots ,a_n\}$ be a generating set of ${\cal U},$
$\psi :\cal U \rightarrow \cal V$ an isomorphism, 
$b_i=\psi (a_i), i=1,\ldots ,n,$ and $B=\{b_1,\ldots ,b_n\}.$
Then $a_i^t=b_i$ is a defining relation in $\cal H.$ 
 Let $c$= $\max \{\left |a_1\right |,\ldots
 \left |a_n\right |, \left |b_1\right |,\ldots \left |b_n\right
|\}.$ 

One can consider two metrics on ${\cal U}$, one is induced from ${\cal G}$ and
another is the word metric in the generators from $A$ (for $W\in {\cal U}$ 
we denote the latter by $\left | W\right | _{\cal U}$.)
We also denote the length of the word $W\in {\cal U}$ 
in the generators $a_1,\ldots a_n$
by ${\left |\left | W \right |\right |}_{\cal U}.$ 
Then   
${\left | W\right |}_{\cal U}
\leq {\left |\left | W \right |\right |}_{\cal U}\leq \left |\left | W 
\right |\right | 
\leq c {\left |\left | W \right |\right |}_{\cal U}.$
It is also clear that ${\left | W\right |}_{\cal U}\leq
\lambda {\left | W\right |},\ 
{\left | W\right |}_{\cal V}\leq\lambda {\left | W\right |}$ and 
${\left | W\right |}\leq c{\left | W\right |}_{\cal U}.$
For $V\in {\cal V}$ we have  
 ${\left | V\right |}\leq c {\left | V\right |}_{\cal V}.$
If $U\in \cal U$ and $U(a_1(J),\ldots ,a_n(J))$ is a geodesic word
in $\cal U$ in the generators in $A,$ then we say that $U$ is
$\cal U$-{\em geodesic} word ({\em or geodesic in} $\cal U$). The same for $\cal 
V.$

If we have a path $p$ in the Cayley graph of ${\cal G}$ such that $\phi (p)$
is a word in the generators in $A$ (resp. $B$) which are words in the
generators of $G$,
then the vertices of $p$ corresponding
to the beginnings and ends of the generators in $A$ (resp. $B$) 
will be called phase vertices (phase points). Phase vertices depend on the way
we write $\phi (p)$ as a word in the generators of $\cal U$ (respectively
$\cal V$).
We call a word $W$ {\em cyclically minimal} if it has minimal length
among all words conjugated to $W$ in ${\cal G}.$

Let $H_1$ be the constant obtained for the group $\cal G$ as in Lemma
\ref{1.9}.
\begin{lm}
\label{Ol2} Let $H_2=H_1+c,\ H=2H_2+8\delta $.  There is a constant
$M_0=M_0({\cal G},{\cal U},{\cal V})$ such that
for any  ${\cal U}$-geodesic words $u, \bar u,$  
${\cal V}$-geodesic words $ v, \bar v,$
and geodesic words $X, Y\in {\cal G},$ \begin{description}
\item[1)] the equality 
$XuY\bar u=1$ for $\left | u\right |\geq M_0$ implies
$$ 
4\cdot \max(\left | X\right |, \left | Y \right |)\geq\left | u\right |, \left 
|\bar u\right |,$$
or $X$ and $Y$ belong to ${\cal U}.$
\item[2)] the equality 
$XuYv=1$ implies $\max (\left | u\right |,\left | v\right |)< M_0$ 
or $ 
4\cdot \max(\left | X\right |, \left | Y \right |)\geq\left | u\right |, \left 
|v\right |,$
\item[3)] the equality 
$XvY\bar v=1$ for $\left | v\right |\geq M_0$ implies one of the following
\begin{enumerate}
\item
$
4\cdot \max(\left | X\right |, \left | Y \right |)\geq\left | v\right |, \left 
|\bar v\right |,$ \item
 $X$ and $Y$ belong to ${\cal V},,$ \item  there are elements $T_1$ and $T_2$
such that $\left | T_i\right | <H ,$
$v\equiv v_1v_3v_2, \bar v\equiv \bar
v_2\bar v_3\bar v_1,$
$$Xv_1T_1^{-1}\bar v_1=1,\ T_1v_3T_2^{-1}\bar v_3=1,\ 
T_2v_2Y\bar v_2=1$$
and
$\left | v_1\right | ,\left |\bar v_1\right | \leq 4 \max\{\left | X\right | 
,H\},
\left | v_2\right | ,\left |\bar v_2\right | \leq 4 \max\{\left | Y\right | 
,H\}$
\end{enumerate}
\end{description}
\end{lm}

{\em Proof:}  See Fig.1. Let $W_1, W_2\in\{u,\bar u, v, \bar v\},$ 
and suppose we have an equality $XW_1YW_2=1.$ If $X=Y=1$ there
is nothing to prove. Suppose $X\not = 1.$ By the condition
of the lemma, there is a quadrangle $p^1q^1p^2q^2$
in the Cayley graph $\Gamma ({\cal G})$ such that $\phi (p^1)=X,
\phi (q^1)=W_1, \phi (p^2)=Y, \phi (q^2)=W_2,$
the paths $q^1, q^2$ are quasigeodesic and $p^1, p^2$ are
geodesic. Consider geodesic paths $s^i,$ such that
$$s^i_{\pm }=q^i_{\pm }.$$ If either $\left | W_1\right |$ or 
$\left | W_2\right |$ are larger than $4 \max\{\left | X\right | ,
\left | Y\right |\}$ then by Lemma \ref{1.7} the paths $s^1$ and $s^2$
contain $8\delta $-bound subsegments $t^1$ and $t^2$
such that $\left | t^1\right | ,\left | t^2\right | >1/3\left | s^1\right | $
(we suppose that $M_0>8\delta/(7/20-1/3)$). We take these subsegments to be
maximal $8\delta $-bound subsegments (this means that there are no
$8\delta $-bound subsegments $\bar t^1, \bar t^2,$ such that
$\bar t^1$ contains $t^1$, $\bar t^2$ contains $t^2$ and one of these
inclusions is proper). 

 Lemma \ref{1.9} allows us to find
phase points  $o_{1j}, (j=1,2)$ on $q^1,$ and $o_{2j}, (j=1,2)$ on $q^2,$ 
such that 
$$\left | o_{11}-t^1_{-}\right |, \left | o_{12}-t^1_{+}\right |, \left | 
o_{21}-t^2_{-}\right |,
\left | o_{22}-t^2_{+}\right | <H_2.$$ This shows that in the case $W_1=v, 
W_2=\bar v$
(the third case) we have either the first possibility or the third  (we will 
show 
that
the second possibility is a particular case of the third).
The subpath $z^1=o_{11}-o_{12}$ has length
$$\left | z^1\right | >1/3 \left | q^1\right | ,$$
and by Lemma \ref{1.5} every vertex of $z^1$ can be connected with a vertex of 
the path $z^2=o_{21}-o_{22}$ by some path $t$ of length $<4H_2+16\delta .$

Consider now vertices $a_1, a_2, a_3\ldots $ (called phase vertices) 
of the path 
$z^1,$ such that the labels of the 
subpaths $a_i-a_{i+1}$ are graphically equal to generators
of ${\cal U}$ (${\cal V}$) . Similarly, choose phase vertices 
$b_1, b_2, b_3,\ldots $ on $z^2.$ As was noticed above, each $a_i$
can be connected with some vertex $b_k$ by a path $t_i$ of 
length $<2H$ . Therefore the number of different 
labels $\phi (t_i)$ of paths $t_i$ is less then $2\left | J\right | ^{2H}.$
So, for  sufficiently large $M_0$ there exist vertices $a_i$ and $a_j$ such that
$\phi (t_i)\equiv \phi (t_j)=T.$ 

Let $z_1$ be the subpath of $z^1$ connecting
$a_i$ with $a_j,$ and $z_2$ be the subpath of $z^2$ connecting
$b_k$ with $b_l$. Observe that $z_1$ (resp. $z_2$) can be made arbitrary
long if one takes long $z^1$ (resp. $z^2$) .

The label of the closed path $t_i^{-1}(a_i-a_j)t_j(b_l-b_k)$ is the word
$T^{-1}\phi (z_1)T(\phi (z_2))^{-1}$.
If $W_1, W_2\in {\cal U}$ then $\phi (z_1), \phi (z_2)
\in {\cal U}$ and if $q^1$ (and hence $s^1$) is long 
enough (increase $M_0$ if necessary), 
then because ${\cal U}$ is conjugate separated, one can find a suitable $T$ 
 in ${\cal U}.$
Hence $X, Y\in {\cal U}.$ 

If $W_1\in {\cal U}$ and $W_2\in {\cal V}$ then, because $\cal H$ is a 
separated HNN-extension, such a $T$ cannot exist
for long $q_1$ or $q_2$ . $\Box $

A direct consequence of this lemma is the following
\begin{cy}
\label{Ol2cy} Let ${\cal U}$ and ${\cal V}$ be subgroups of ${\cal G}$ as in
Theorem 1.  If
$M>\max \{M_0,4, c\lambda \}$ then
for any ${\cal U}$ -geodesic words $u, \bar u,$ any $\cal V$-geodesic
words $v, \bar v,$ 
and geodesic words $X, Y\in {\cal G},$ \begin{itemize}
\item the equality 
$XuY\bar u=1$  implies
$$
M\cdot \max(\left | X\right |, \left | Y \right |)>\left | u\right |, \left 
|\bar 
u\right |,$$
or $X$ and $Y$ belong to ${\cal U}.$
\item the equality 
$XuYv=1$ implies 
$$\max (\left | u\right |,\left | v\right |)< 
M\ \max(\left | X\right |, \left | Y \right |, 1),$$ 
\item the equality 
$XvY\bar v=1$  implies one of the following
\begin{enumerate}
\item
$
M\cdot \max(\left | X\right |, \left | Y \right |)>\left | v\right |, \left 
|\bar 
v\right |,$ \item
 $X$ and $Y$ belong to ${\cal V}$ \item  there are elements $T_1$ and $T_2$
such that $\left | T_i\right | <H,$
$v\equiv v_1v_3v_2, \bar v\equiv \bar
v_2\bar v_3\bar v_1,$
$$Xv_1T_1^{-1}\bar v_1=1,\ T_1v_3T_2^{-1}\bar v_3=1,\ 
T_2v_2Y\bar v_2=1$$
and
$\left | v_1\right | ,\left |\bar v_1\right | \leq 4 \max\{\left | X\right | 
,H\},
\left | v_2\right | ,\left |\bar v_2\right | \leq 4 \max\{\left | Y\right | 
,H\}$
\end{enumerate}
\end{itemize}
\end{cy}

\begin{df} Suppose a word $W$ in $\cal G$ has the following decomposition
$$W\equiv X_0W_1X_1W_2X_2\ldots W_kX_k,$$
where each $X_i$ is a reduced word, 
each $W_i$ either belongs to ${\cal U}$ and is $\cal U$-geodesic 
 or belongs to ${\cal V}$ and is $\cal V$-geodesic, 
and if, for some $i,$ $W_i$ and $W_{i+1}$ both belong to  ${\cal U}$
or both belong to ${\cal V}$ 
then $X_i\not =1$. Then this decomposition will be called a $\cal U
\cal V$-{\em decomposition of} $W.$

If $W=1,$ $X_0\equiv 1$ 
and in the above definition indices are taken modulo $k,$
then the above decomposition will be called a cyclic $\cal U
\cal V$-{\em decomposition of} $W.$

In the case where all the $W_i$'s belong to $\cal U$ (resp. $\cal V$), we
will talk about $\cal U$-{\em decomposition} ($\cal V$-{\em decomposition}).
\end{df}

\begin{df} A $\cal U\cal V$-decomposition of the word $W$
$$W\equiv W_1X_1W_2X_2\ldots W_kX_k=1$$
is called {\em splittable} if one of the following holds.
\begin {enumerate}
 
\item There is a $j\leq k$ such that 
$W_1\equiv W_{11}W_{12},\  W_{11},W_{12}\in {\cal U} $
and $W_j\equiv W_{j1}W_{j2}, \ W_{j1},W_{j2}\in {\cal U} $
 and 
$$W_{12}X_1W_2X_2\ldots W_{j1}=W'\in {\cal U},$$
where $$\left | W_{12}\right | _{\cal U} +\left | W_{j1}\right | _{\cal U}>
\left | W'\right | _{\cal U}+\left | W_{j1}W'^{-1}W_{12}
\right |  _{\cal U}.$$   

\item There is a $j\leq k$ such that $W_1\equiv W_{11}W_{12}, W_{11},W_{12}\in 
{\cal V} $
and $W_j\equiv W_{j1}W_{j2}, W_{j1},W_{j2}\in {\cal V} $
and 
$$W_{12}X_1W_2X_2\ldots W_{j1}=W'\in {\cal V},$$
where $$\left | W_{12}\right | _{\cal V}+\left | W_{j1}\right | _{\cal V}
>\left | W'\right | _{\cal V}+\left | W_{j1}W'^{-1}W_{12}
\right |  _{\cal V}.$$   

\item There is a $j\leq k$ such that $W_1\equiv W_{11}W_{12}, W_{11}\not =1,
W_{12}\not =1\in {\cal V} $
and $W_j\equiv W_{j1}W_{j2}, W_{j1}\not =1, W_{j2}\not = 1 \in {\cal V} $
and $$W_{12}X_1W_2X_2\ldots W_{j1}=T,$$
where $T\not\in {\cal V},$  $\left | T\right | <H$
and $W_{12}, W_{j1}$ have minimal length in $\cal V$ 
among the subwords of $W_1$ and $W_j$ with the property, that the
above equality is satisfied for some word $T,$ where $\left | T\right | <H$.
\end {enumerate}

A cyclic $\cal U\cal V$-decomposition of $W$ 
is called  {\em cyclically nonsplittable} if all  permutations
of the form 
$$W_iX_i\ldots W_kX_kW_{1}\ldots W_{i-1}X_{i-1}=1.$$
are nonsplittable.
\end{df}

Our main goal in this section is to prove Corollary 6.

The following lemma follows from Corollary \ref{Ol2cy}.
\begin{lm}
\label{Ol3} 
For any $\cal U$ -geodesic words $u, \bar u\in {\cal U},$ and
$\cal V$-geodesic words $v, \bar v\in {\cal V}$ 
and geodesic words $X, Y\in {\cal G}$ we have the following :
\begin{itemize}
\item 
If $XuY\bar u=1$ is a cyclically nonsplittable $\cal U$-decomposition,
then
$$
M\cdot (\left | X\right | + \left | Y \right |)>\left | u\right |, \left |\bar 
u\right |.$$
\item The equality 
$XuYv=1$ implies $$\max (\left | u\right |,\left | v\right |)< 
M\ \max(\left | X\right | + \left | Y \right |, 1),$$ 
\item  
If $XvY\bar v=1$ is a cyclically nonsplittable $\cal V$-decomposition
then 
$$
M\cdot (\left | X\right | + \left | Y \right |)>\left | v\right |, \left |\bar 
v\right |.$$ 
\end{itemize}
\end{lm}

\begin{prop}
\label{Pr}
Suppose  that in the group ${\cal G}$ 
we have a cyclically nonsplittable $\cal U\cal V$-decomposition
of the word $W:$
$$W\equiv W_1X_1W_2X_2\ldots W_kX_k=1,$$
where the $X_i$ are geodesic words.

Then  for any $i$ $$\left | W_{i}\right |
\leq M(\sum_{i=1}^k\left | X_i\right |)+2M(k-1)(2H_2+6\delta)+k(H_2+1) .$$
(See Fig. 3.)
\end{prop}

We will give a proof of 
this proposition together with the following lemma
by simultaneous induction on $k.$
\begin{lm} Suppose we are given a $\cal U\cal V$-decomposition of 
a word $\bar W:$
$$ \bar W\equiv\bar X_0\bar W_1\bar X_1\bar W_2\bar X_2\ldots \bar W_{k-2}\bar 
X_{k-2}\bar W_{k-1},$$
$k>2$, and $\bar X_0=1$ if $k>3,$ 
the $\bar X_i$'s are geodesic, 
each $\bar W_i$ 
either belongs to ${\cal U}$ and is geodesic in $\cal
U,$ or
belongs to $\cal V$ and is geodesic in $\cal V$. 

Let 
$U_k$ be a geodesic word such that $U_k=\bar W$
and $\bar U_k$ be a geodesic such that
$$\bar U_k=U_k\bar W_{k-1}^{-1}.$$
Suppose that 
the decomposition
$$\bar X_0\bar W_1\bar X_1\bar W_2\bar X_2\ldots \bar W_{k-2}\bar X_{k-2}
\bar W_{k-1}U_k^{-1}=1,$$
is cyclically nonsplittable. 
Let  
$\bar V_i$ be a geodesic word such that $\bar V_i=\bar W_i.$
Then Gromov's products $$L_k=(\bar U_k^{-1}
\cdot \bar V_{k-1})\leq
M(\sum_{i=1}^{k-2}\left | \bar X_i\right 
|)+M(1+2(k-2))(2H_2+6\delta)+(k-1)(H_2+1) 
.$$
(See Fig. 2, 4.)
\end{lm}
The assertion of the Proposition for $k=2$ is just Lemma \ref{Ol3}.

We first will prove Lemma 6 for $k=3.$
The proof is illustrated in Fig.2.
Consider a geodesic triangle $[x_1,x_2,x_3]$ in $\Gamma ({\cal G})$ such
that the label $\bar A$ of the side $[x_1,x_2]$ is equal to $\bar U_3^{-1}$
in ${\cal G}$ and the label $\bar B$ of the side $[x_1,x_3]$ is equal
to $\bar V_{2}.$
Then $L_3=(\bar A\cdot \bar B).$
By Lemma~\ref{1.5}, one can find decompositions $\bar A\equiv A_1A_2,\ \  
\bar B\equiv B_1B_2,$ where $\left | A_1\right | =\left | B_1\right | >L_3-1,$ 
and
for subsegments $p,q$ of the sides $[x_1,x_2]$ and $[x_1,x_3]$ with
labels $A_1$ and $B_1,\ \ \left | p_+-q_+\right |\leq 4\delta .$
Consider the path $s$ issuing from $x_1$ and having $\bar W_2$ as a label.
According to Lemma \ref{1.9} there is an initial subpath $s^1$ of $s$
with the label $\bar W_{21}$ belonging to the same subgroup $\cal U$ or
$\cal V$ as $W_2$, such that $\left | q_+-s^1_+\right | <H_1+c=H_2.$ Let
$t$ be a path issuing from $x_1$ and having the label $\bar X_1^{-1}\bar  
W_1^{-1}.$
Then either $\left | t_+-p_+\right | <\delta +\left | \bar X_0\right | $ or 
there 
is a subpath
$t^1$ issuing from $x_1$ and having the label  $\bar X_1^{-1}\bar W_{12}^{-1},$
such that $\left | t^1_+-p_+\right | <2\delta +H_2.$

It is clear that $L_3-1\leq H_2+\left | \bar W_{21}\right |.$

By Corollary \ref{Ol2cy}, $\left | \bar W_{21}\right | <M \max\{\left | \bar 
X_1\right |, 
H_2+5\delta +\left | \bar X_0\right |, 2H_2+6\delta, H_2+6\delta +\left | \bar 
X_1\right | , 1\}.$
Hence $L_3\leq H_2+1 + M((\left | \bar X_0\right | +\left | \bar X_1\right | 
)+2H_2+6\delta).$
The Lemma is proved for $k=3.$

Suppose now that  the Proposition is
proved for all $k< r$ and the Lemma is
proved for all $k\leq r.$  We will prove the Proposition for $k=r.$ 
Let $V_{i}$ be geodesic words such that $V_{i}= W_i.$
See Fig. 3.
Let $[x_1,y_1,\ldots x_r,y_r]$ be a geodesic $2r$-gon in the 
Cayley graph $\Gamma ({\cal G}),$ such that $V_{i}$ is the label of $[x_i,y_i]$
and $X_i$ is the label of $[y_i,x_{i+1}]$ 
($X_r$ is the label of $[y_r,x_1]$). Let $U$ be the label of $[y_1,x_3],$
and $W$ the label of $[x_1,x_3].$ 

By the assertion of the Lemma for $k=r,$ applied to the word
$$W_3X_3\ldots W_rX_rW_1,$$
we have 
$$L_{r}=(W\cdot V_{1})\leq (r-1)(H_2+1)+ M(\sum_{i=3}^r\left | X_i\right | +
(1+2(r-2))(2H_2+6\delta )).$$ 

By the assertion of the Lemma for $k=3,$ applied to the word
$$X_2^{-1}W_2^{-1}X_1^{-1}W_1^{-1},$$
we have $$L_3=(U\cdot V_{1}^{-1})\leq H_2+1 + M((\left | X_0\right | +\left | 
X_1\right | 
)+2H_2+6\delta).$$
Now, by Lemma~\ref{1.5}, $$\left | V_{1}\right | =
 L_3+L_r\leq r(H_2+1)+ M(\sum_{i=1}^r\left | X_i\right | +(r-1)2(2H_2+6\delta 
)). 
$$
 
 Now we will prove the Lemma for $k=r+1.$ See Fig. 4.
 
 Consider a $2r$-gon $[x_1,y_1,\ldots ,x_{r},y_{r}]$ in $\Gamma ({\cal G})$
such that $\bar X_i$ is the label of $[y_i,x_{i+1}]$ and
a path with the label $\bar W_i$ has initial point $x_i$ and terminal
point $y_i.$ 
 
 Consider a geodesic triangle $[x_1,x_{r},y_{r}]$  such
that the label $\bar A$ of the side $[x_{r},x_1]$ is equal to $\bar U_r^{-1}$
in ${\cal G},$ and the label $\bar B$ of the side $[x_{r},y_{r}]$ is
equal to $\bar V_{r}.$
Then $L_r=(\bar A\cdot \bar B).$
By Lemma \ref{1.5}, one can find decompositions $\bar A\equiv A_1A_2,\ \  
\bar B\equiv B_1B_2,$ where $\left | A_1\right | =\left | B_1\right | >L_r-1$ 
and
for subsegments $p,q$ of the sides $[x_{r},x_1]$ and $[x_{r},y_{r}]$ with
labels $A_1$ and $B_1,\ \ \left | p_+-q_+\right |\leq 4\delta .$
Consider the path $s$ issuing from $x_{r}$ and having $\bar W_{r}$ as its label.
According to Lemma \ref{1.9} there is an initial subpath $s^1$ of $s$
with the label $\bar W_{r1}$ belonging to the same subgroup $\cal U$ or 
$\cal V$ as $W_r,$ such that $\left | q_+-s^1_+\right | <H_2.$ 
There are the following $r$ possibilities:
$\left |p_+,[x_1,y_1]\right |<\delta ,$
$\left |p_+,x_{i+1}\right |<i\delta +\left | \bar X_{i}\right |,$
($i=1,\ldots r-1$).

We consider the first possibility. Then there 
is a path $t$ issuing from $y_1$ with the label $\bar W_{12}^{-1}$
such that $\left | p_+,t_+\right | <\delta +H_2.$ Let $f$ be the label of the 
geodesic path $s^1_+,t_+,$ then $\left | f\right | <5\delta +2H_2.$

 If $f\not =1$ or both $\bar W_{12},\bar W_{r1}\not\in {\cal U} ({\cal V}),$ 
then
the word  $$\bar W_{12}\bar X_1\ldots \bar W_{r-1}\bar X_{r-1}\bar W_{r1}f=1$$ 
is
nonsplittable and  
 we can apply the Proposition for $k\leq r$ to get
 $$L_{r+1}-1\leq r(H_2+1)+ M(\sum_{i=1}^{r-1}\left | \bar X_i\right | 
+1+2(r-1)(2H_2+6\delta )). $$
 
 Now suppose that $f=1, \bar W_{12},\bar W_{r1}\in {\cal U}.$
Let $\hat W$ be a quasigeodesic word such that
$\hat W=\bar W_{r1}\bar W_{12}.$
 Consider instead of the word
 $$\bar W_{12}\bar X_1\ldots \bar W_{r-1}\bar X_{r-1}\bar W_{r1}=1$$
 the word
$$\hat W\bar X_1\ldots \bar W_{r-1}\bar X_{r-1}$$ 
(it is nonsplittable) ,
and apply the Proposition for $k=r-1$ to estimate $\hat W$

Because nonsplitability implies
${\left | (\bar W_{r1}\bar W_{12})\right |}_{\cal U} = \left | \bar W_{r1}\right 
| 
_{\cal U}+
\left | \bar W_{12}\right | _{\cal U},$ and hence 
${\left | \bar W_{r1}\right |}\leq {\left | \hat W\right |}+H_2,$
the proof is finished in this case.  

The case $f=1, \bar W_{12},\bar W_{r1}\in {\cal V}$ can be considered
similarly.

The other $r-1$ possibilities $\left |p_+,x_{i+1}\right |<i\delta +\left | \bar 
X_{i}\right |,$
($i=1,\ldots r-1$) can be considered similarly to the first possibility.
Proposition 1 and Lemma 6 are proven. $\Box $  
\begin{cy}
\label{est} There are constants $M_1$, $M_2$ such that if
in the group ${\cal G}$ 
we have a cyclically nonsplittable 
$\cal U\cal V$-decomposition
$$W\equiv W_1X_1W_2X_2\ldots W_kX_k=1,$$
then for $W_i\in\cal U$ we have
${\left | W_i\right | }_{\cal U}\leq M_1\sum_i \left | X_{i}\right | + M_2k$
and for  $W_i\in\cal V$ we have
${\left | W_i\right | }_{\cal V}\leq M_1\sum_i \left | X_{i}\right | + M_2k.$
\end{cy}
\section{Diagrams}
Recall that a map a is finite,
planar connected $2$-complex. 

By a {\em diagram} $\Delta$ over a presentation 
$<a_1,\ldots ,a_m |
R_1,\ldots ,R_n>,$ where the words $R_i$ are cyclically reduced,
we mean a map  with a function $\phi$ 
which assigns to each edge of the map one of the letters $a_k^{\pm 1},$
$1\leq k\leq m,$ such that $\phi (e^{-1})= ({\phi (e)})^{-1}$ and if
$p=e_1\ldots e_d$ is the contour of some cell $\Phi$ of $\Delta,$
then $\phi (p)\equiv\phi (e_1)\ldots \phi (e_d)\equiv R$ is a cyclic
shift of one of the defining words $R_i^{\pm 1}.$ In general the word
$\phi (p)$ is called the {\em label} of the path $p.$ The label of a diagram
$\Delta$ (whose contour is always taken with a counterclockwise
orientation) is defined analogously.

Van Kampen's lemma states that a word $W$ represents the identity
of the group $G$ if and only if there is a simply   connected
(or {\em Van-Kampen}, or {\em disk}) diagram $\Delta$ over $G$ such that the 
boundary
label of $\Delta$ is $W.$


Due to Van-Kampen's lemma, a group $G$ is hyperbolic if and only if there  are 
constants $K$ and $C$ such that for any element $W=1$ in ${\cal G}$ there
is a diagram with  boundary label $W$ and number of cells
$\leq K\left |\left | W\right |\right | +C.$ Since  $\left |\left | W\right 
|\right |\geq 1$
we can assume (taking $K+C$ instead of $K$) that $C=0.$  

\section{HNN-extensions}
Let ${\cal H}$ be the HNN-extension as in Theorem 1. 
We fixed the presentation for the group $\cal H$ in the first section.
As we just noticed, from the hyperbolicity  of the group ${\cal G}$  
it follows that there is a 
constant $K$ such that for any element $W=1$ in ${\cal G}$ there
is a diagram with boundary label $W$ and number of cells
$\leq K\left |\left | W\right |\right | .$ 

The contents of this section is the proof of the following 
\begin{prop} There is a linear function $L(x)$ depending only on ${\cal G},
{\cal U},{\cal V}$ such that
for any element $W=1$ in ${\cal H}$ there
is a diagram over ${\cal H}$ with boundary label $W$ and number of cells
$\leq L(\left |\left | W\right |\right | ).$
\end{prop}

Let $\Delta $ be a diagram over ${\cal H}$ with
boundary label $W$. New cells 
corresponding to the relations $t^{-1}a_it=b_i,$ where 
$i=1,\ldots ,n,$
will be called $t$-cells.
They are shown on Fig. 5. A 
configuration of $t$-cells in a
diagram over ${\cal H},$ as shown on Fig. 6a, we call a {\em $t$-annulus}. 

From now on we suppose that $\Delta$ is {\em minimal}, this means that it 
has a minimal possible number of $t$-cells.

\begin{lm} 
A minimal diagram over ${\cal H}$ cannot contain a $t$-annulus.
\end{lm}
{\em Proof:} Suppose it contains a $t$-annulus. Take a $t$-annulus
such that there is not another $t$-annulus inside it. Then the 
label of its internal contour $p$  equals the identity in the 
group ${\cal G}.$ Hence 
the 
label of its external contour $q$  equals the identity in the 
group ${\cal G}.$
We can decrease the number of $t$-cells by applying the following

{\bf Transformation  1} Assuming that the contour $p$ in Fig. 6a
bounds a $\cal G$-diagram, replace the interior of the diagram having
the  contour $q$ in Fig. 6a  
by a $\cal G$-diagram with the contour $q.$

The Lemma is proved.

By this lemma, 
$t$-cells can only  form
 $t$-strips as shown in Fig. 6b, and these
$t$-strips must end on the boundary of $\Delta .$

The direction of the $t$-edges defines $\cal U$-  and $\cal V$-sides of a 
$t$-strip.
The minimality of $\Delta$ also implies that the $\cal U$- ($\cal V$-)sides 
of the 
$t$-strips are geodesic words
respectively in the subgroups $\cal U$ and $\cal V$. 
Indeed, suppose we have a $t$-strip, such that the 
$\cal U$- and $\cal V$-sides 
of this
$t$-strip are not geodesic words
in the subgroups $\cal U$ and $\cal V$. Let the path $p$ correspond
to the $\cal V$-side of the $t$-strip, with $\phi (p)=V_1.$ Let $V_2$ be
a geodesic word in $\cal V$ such that $V_2=V_1$ in $\cal V.$
Let $V_2=\phi (q).$ We replace the $t$-strip by the diagram $\Theta $
having the same contour as shown in Fig. 7 and cut out the annulus.
As a result we have a diagram with fewer $t$-cells.

The typical form of $\Delta$ is shown in Fig. 8. 

Our diagram $\Delta$ is subdivided by the $t$-strips into a  set
of mutually disjoint maximal ${\cal G}$-subdiagrams.
The maximal ${\cal G}$-subdiagrams are the connected components of $\Delta$
which remain after deleting all $t$-edges and interiors of $t$-cells.

Our next goal is to study maximal ${\cal G}$-subdiagrams in the diagram 
${\Delta}$.
 A typical form of a maximal ${\cal G}$-subdiagram is shown in Fig.9.

\begin{df}
An {\em island} in a maximal ${\cal G}$-subdiagram is a 
${\cal G}$-subdiagram with the following properties:
\begin{enumerate} 
\item The contour of the   island is subdivided into paths, such that
every path either belongs to the boundary of $\Delta$ or is a part
of a ${\cal U}$- or ${\cal V}$-side   of a $t$-strip in $\Delta.$
\item There is no point on the contour of the   island such that the deletion of 
this point splits the   island into two or more disconnected components.
\end{enumerate}
\end{df} 
Every edge of an  island is {\em proper}, i.e. lies on the boundary 
of 
some 
cell of this island.

A {\em bridge} is a  pair of paths $\{p,p^{-1}\},$ where  $p=e_1\ldots e_r$ 
is a maximal subpath, consisting of improper edges of $\Delta$, such that the 
valencies of the terminal points of $e_1\ldots e_{r-1}$ are equal to 2.

Every maximal $\cal G$-subdiagram consists of islands and bridges.

The contour of each island is canonically
subdivided into paths, and each of these paths either
belongs to the ${\cal U}$- or ${\cal V}$-side of some $t$-strip 
and is maximal with respect to this property (we will call these paths
{\em ${\cal U}$-paths} and {\em ${\cal V}$-paths}, or generally {\em strippaths} 
or belongs to the boundary of
$\Delta$ and is maximal with respect to this property (call them {\em boundary
paths}.)

We will call a vertex on a ${\cal U}$-side (resp. ${\cal V}$-side) 
of a $t$-strip a {\em 
phase vertex} if
it corresponds to the beginning or end of a path labelled by 
some generator $a_i$ of ${\cal U}$ 
(resp. generator $b_i$ of ${\cal V}$) and to the beginning or to the end  
of the $t$-edge.  A vertex on a strippath will be called 
 a {\em 
phase vertex} if it is a phase vertex on the corresponding $\cal U$- or 
$\cal V$-side of the $t$-strip.

{\bf Transformation 2}
We now carry out  surgeries on the diagram. 
Our objective is
 to make ${\cal U}$- and ${\cal V}$-paths contact the boundary paths only 
through
phase vertices.  
 See Fig. 10. 
Let $p$ be a boundary-path of an island $I$ 
which is adjacent to 
a ${\cal U}$-path at the point $O_1$. Let $O_2$ be the 
phase vertex on the ${\cal U}$-path closest to $O_1$ 
and $q$ be a subpath of the ${\cal
U}$-path connecting $O_1$ and $O_2.$ We make a cut along the path $q.$ 
Each cut can increase the boundary by most $2c.$ Collectively, these cuts define 
Transformation 2. 

Such a transformation increases the length of the boundary of $\Delta $
by a factor of not more then $(1+4c).$
   This coefficient
does not depend on the diagram. So, if we
can prove a linear isoperimetric unequality for a transformed diagram (or
diagrams), we
can prove it for the original one. 

Without loss of generality 
we assume now that $\Delta$ has the property that all boundary paths
contact strip-paths only at phase vertices.

{\bf Transformation 3}
For each boundary-path $q$  on the contour of an island, linking
two ${\cal V}$-paths and 
such that $\phi (q)$ is equal in ${\cal G}$ to an
element of ${\cal V}$, $\phi (q) =v,
$
as shown on Fig. 11, glue to $q$ two diagrams over ${\cal G}$: 
${\Theta}_{q,1}$ with contour $q^{-1}p$ and ${\Theta}_{q,2}$ with  contour
$p^{-1}q$, where $\phi (p)\equiv v$ and $v$ is a geodesic word in ${\cal V}$.
  See Fig. 12.

 We do the same for each  boundary-path $q$   such that $\phi (q)$ is 
equal in ${\cal G}$ to an element of ${\cal U}$ and 
linking    
two ${\cal U}$-paths.
The resulting diagram will have the same contour as $\Delta.$
Transformation 3 ends by cutting out each diagram ${\Theta}_{q,2}.$ 

Since $\cal U$ and $\cal V$ are quasiisometrically embedded 
in ${\cal
G}$ and ${\cal
G}$ is
hyperbolic, we can pick each 
diagram ${\Theta}_{q2}$ over ${\cal G}$, 
with a contour $pq^{-1},$ where $\phi (p)=v (u)$,  
so that the number of cells in it is 
less than $K(\left |\left | v\right |\right | + 
\left |\left | q\right |\right | )\leq K(c{\lambda}+1)\left |\left | q\right 
|\right |
.$
Hence the sum of cells in all the diagrams ${\Theta }_{q,2}$ is less than 
$K(c{\lambda}+1)\left |\left | W\right |\right |.$

  Our goal is to bound the number of cells in the diagrams  with contour
$W$ by a linear function of $\left |\left | W\right |\right |.$ 
We have bounded the number of cells in 
the union of the diagrams ${\Theta }_{q,2}$ for all boundary paths $q$
of the type considered above. 
The length   of the contour of the resulting 
diagram is less then $c{\lambda}$ times the length of the contour of the  
original diagram
$\Delta $. 
Now instead of the diagram $\Delta$  
we will consider this new diagram which we will also denote by 
${\Delta}.$

From now on we do not change the boundary of $\Delta$ anymore.

\begin{lm} 
\label{phase}
If 
in the diagram $\Delta$ on the boundary of an island there are
two ${\cal U}$-
(${\cal V}$)-paths $p$ and $s$ such that the terminal vertex of $p$
is the initial vertex of $s$ and is a phase vertex for both $p$ and $s$,  
then the path $ps$ 
is a geodesic in the group ${\cal U}$
(${\cal V}$).
\end{lm} 

{\em Proof:} See Fig. 13. Suppose that on 
the boundary of an island two ${\cal U}$-paths 
$p$ and $s$ have
a common phase vertex $p_+ =s_-$ 
with $\phi (p)=u_1$, $\phi (s)=u_2.$
Suppose that their union is not a
geodesic path in ${\cal U},$ and let $q$ be a geodesic path in $\cal U$
such that 
 $\phi (q)=u$ and $u=u_1u_2.$
So ${\left | q\right | }_{\cal U} < {\left | p\right | }_{\cal U}+{\left | 
s\right 
| }_{\cal U}.$

We make the 

{\bf Transformation 4} as shown on Fig. 13 (the contour of
the subdiagram is not changed by this transformation). 
We cut along the path $ps$ and incert two mirror copies of the
diagram with the contour $psq^{-1}$. Then
we cut along the
edge $q$ and insert a patch of two adjoining $t$-strips.
We then cut 6 $t$-cells and reattach them in a different way to create 
the subdiagram in Fig. 13c. 
Transformation 4 ends by cutting out the $t$-annulus (Transformation 1).

After the cutting we have a diagram with fewer $t$-cells. This
contradicts the minimality of $\Delta.$ Indeed, we replace 
${\left | p\right | }_{\cal U}+{\left | s\right | }_{\cal U}$ $t$-cells by 
${\left | q\right | }_{\cal U}$ $t$-cells.
This completes the proof of the Lemma.

In the diagram ${\Delta},$
${\cal U}$-sides of two t-strips 
cannot be glued together along a path longer than $Mc+c$. If they are glued from
one common phase vertex to another one then we can make a $t$-annulus
then cut it out using Transformation 1 and decrease the number of $t$-cells. 
This contradicts the minimality of $\Delta .$
If they are glued not from
one common phase vertex to another one then we can apply Corollary \ref{Ol2cy}
and restrict their length by $Mc+c$.

\begin{df}
A $\cal G$-subdiagram 
with a cyclic $\cal U\cal V$-decomposition of boundary label 
$$W\equiv W_1X_1W_2X_2\ldots W_kX_k=1,$$
where 
each $W_i$ is a label of a strippath,  is called nonsplittable if the
above decomposition is
cyclically nonsplittable (see Definition 3).
\end{df}

\begin{df}
A maximal nonsplittable $\cal G$-subdiagram is called a nonsplittable piece.
\end{df}

\begin{df} A $\cal V$-piece 
is a $\cal G$-subdiagram having 
boundary label 
$$W_1X_1W_2X_2=1,$$ where $W_1$, $W_2$ label only subpaths of  
${\cal V}$-sides of $t$-strips, 
and $X_1, X_2$ are  shorter than $H$ (the $X_i$'s can be trivial),
and we assume that the $\cal V$-piece is not
properly contained in another $\cal G$-subdiagram
with a boundary label of that type.\end{df}

\begin{df} By a piece we mean either a nonsplittable piece or a 
${\cal V}$-piece.\end{df}

A piece may consist of a several islands and not necessarily be an 
island itself.

A contour of a piece consists of strippaths (which are maximal with respect 
to the property of belonging to a side of a $t$-strip and to the contour of the 
piece, and beginning and ending in the phase vertex) and paths connecting them.
\begin{lm}
\label{na}
Every maximal $\cal G$-subdiagram in $\Delta$ 
is partitioned into nonsplittable pieces and 
${\cal V}$-pieces between them.
\end{lm}

{\em Proof:} Let $\Theta$ be a maximal $\cal G$-subdiagram. A boundary label of 
$\Theta $ is a word $$W\equiv W_1X_1W_2X_2\ldots W_kX_k=1,$$
where 
$W_i\in\{{\cal U},{\cal V}\}$ and is geodesic in ${\cal U}$ or in ${\cal V}$. 
Also if $W_i,W_{i+1}\in {\cal U}$ or $W_i,W_{i+1}\in {\cal V}$ 
(indices are taken modulo $k$)
then $X_i\not=1.$

Suppose the boundary label is a splittable $\cal U\cal V$-decomposition.  
Suppose that the second possibility in Definition 2 holds.
Then there is a $j\leq k$ such that $W_1\equiv W_{11}W_{12}, W_{11},W_{12}\in 
{\cal V} ,$
$W_j\equiv W_{j1}W_{j2}, W_{j1},W_{j2}\in {\cal V} $
and 
$$W_{12}X_1W_2X_2\ldots W_{j1}=W'\in {\cal V},$$
where $$\left | W_{12}\right | +\left | W_{j1}\right | >\left | W'\right | 
+\left 
| W_{j1}W'^{-1}W_{12}
\right |.  $$   
Let $\bar W_{11}, \bar W_{12}, \bar W_{j1}, \bar W_{j2}$ be the words
corresponding to $W_{11}, W_{12}, W_{j1}, W_{j2}$ on the other side
of the $t$-strips.
Then we can make 

{\bf Transformation 5 }  
as shown in Fig. 14. 
Let $V$ be a geodesic word in $\cal V$ such that $V=W_{j1}W'^{-1}W_{12}.$
Let $\bar V$ be the word corresponding to $V$ on the other side of the 
$t$-strip.
Let ${\Delta}_1$ be the subdiagram with the contour $$t^{-1}{\bar W}_1t
V^{-1}t^{-1}{\bar W}_{j}tW_{j2}^{-1}W'^{-1}W_{11}^{-1}.$$
Replace it by the union of three diagrams:
${\Theta}_1,$ which is just 
a $t$-strip with the boundary label $tV^{-1}t^{-1}{\bar V}$,
${\Theta}_2,$
with the boundary label $${\bar V}^{-1}tt^{-1}{\bar W}_{j1}tt^{-1}{\bar W}'^{-1}
tt^{-1}{\bar W}_{12}tt^{-1}$$ and
${\Theta}_3,$ with the
boundary label $${\bar W}_{11}tt^{-1}{\bar W}'tt^{-1}{\bar 
W}_{j2}tW_{j2}^{-1}W'^{-1}
W_{11}^{-1}t^{-1}.$$

${\Theta}_1$ is glued to ${\Theta}_2$
along the path with the label $\bar V$, 
${\Theta}_3$ is glued to ${\Theta}_2$
along the path with the label $\bar W'$. The union of  ${\Theta}_1,$
${\Theta}_2$ and ${\Theta}_3$ has the same boundary label as ${\Delta}_1.$

Now instead of one maximal $\cal G$-subdiagram $\Theta$ we obtained
three: ${\Theta}_4, {\Theta}_5$ and ${\Theta}_6$.
where ${\Theta}_4$ has 
the boundary label $X_1\ldots X_{j-1}V$, ${\Theta}_5$
has the boundary label $V^{-1}W_{j1}W'^{-1}W_{12}$ and ${\Theta}_6$ has the
boundary label $W_{11}W'W_{j2}X_j\ldots X_k$. 
The diagram ${\Theta}_5$ is the interior of a $t$-annulus, and, together with 
this $t$-annulus,  gives ${\Theta}_2$. We can end Transformation~5 by 
application
of Transformation 1 and replacing ${\Theta}_2$ 
by a diagram over $\cal U.$
This decreases the number of $t$-cells in $\Delta ,$  because instead of
${\left | W_{12}\right |}_{\cal V}+{\left | W_{j1}\right |}_{\cal V}$ $t$-cells
we have now ${\left | V\right |}_{\cal V}+{\left | W'\right |}_{\cal V}$.
This contradicts the assumption that $\Delta$ is minimal.

The case where the first possibility in Definition 3 holds, can be considered
similarly.

Suppose now that the third possibility holds. 
Then there is a $j\leq k$ such that $W_1\equiv W_{11}W_{12},$ with $ W_{11}
\not =1$
and $W_{12}\not =1$ both belonging to ${\cal V}, $
$W_j\equiv W_{j1}W_{j2},$ with $ W_{j1}\not =1$ and  $W_{j2}\not = 1$
both belonging  to 
$ {\cal V} ,$
and $$W_{12}X_1W_2X_2\ldots W_{j1}=T,$$
where $T$ is not necessarily in ${\cal V},$ but $\left | T\right | <H$
and $W_{12}, W_{j1}$ have minimal length among the subwords with this property
(see Fig. 15).
Then we can represent $\Theta$ as a union of two 
$\cal G$-subdiagrams ${\Theta}_1$ and  ${\Theta}_2$ where ${\Theta}_1$
has   contour
label  
$$W_{12}X_1W_2X_2\ldots W_{j1}T^{-1},$$ and ${\Theta}_2$ has 
contour label
$$W_{11}TW_{j2}\ldots X_k.$$ 
If $W_{j2}=W_{j21}W_{j22},$ 
$W_{11}=W_{111}W_{112}$ and $W_{j22}\ldots X_kW_{111}=T_1,$ 
where $T_1\not\in {\cal V}$ but $\left | T_1\right | <H,$
then $\Theta$ is the union of ${\Theta}_1$ and a $\cal G$-subdiagram 
${\Theta}_3$
with contour $W_{111}T_1W_{j22}X_j\ldots X_k,$
connected to $\Theta _1$ by a $\cal V$-piece having contour 
$TW_{j21}T_1^{-1}W_{112},$ 
as shown on Fig.15. 

Every time when a diagram under consideration is splittable we represent it as
a union of several subdiagrams. Continuing this process we will obtain
the desired partition.
$\Box$
\begin{df}
A piece  
with a contour that does not contain any boundary paths, is called 
a {\em concealed piece}. 
\end{df}
\begin{df}
A piece is called a {\em $k$-piece} if it
has k strippaths on the contour.\end{df}
\begin{lm} 
\label{NK} Suppose we have fixed some partition of maximal 
${\cal G}$-subdiagrams in $\Delta $ into nonsplittable pieces and 
$\cal V$-pieces.
Let $N_k$ ($k>2$) be the number of nonsplittable $k$-pieces 
and $N_2$ be the number of nonsplittable $2$-pieces 
plus the number of 
$\cal V$-pieces
in
 $\Delta,$  and let $S$ be the number of
$t$-strips in it, $S\geq 3.$ 
Then  $J:=\sum_{k=2}^SN_kk\leq 15 (S-2).$
\end{lm}

{\em Proof:} 
The number of concealed pieces between two $t$-strips
cannot be more than one. 

We say that two $t$-strips {\em adjoin irregularly} if they do not adjoin 
through phase vertices. A subdiagram with a contour formed by two $t$-strips
adjoined irregularly is an island, but it is always either properly included
in some nonsplittable piece or included in some ${\cal V}$-piece.

From now on, we will not worry about filling in the pieces by subdiagrams
over ${\cal G},$ we just consider all possible configurations of
$S$ $t$-strips in the plane. We can forget that the strips are $t$-strips
and think about them simply as about strips.
We treat the paths shorter than $H$ on the boundaries
of  $\cal V$-pieces as if they were just points, 
so that at these points the corresponding $t$-strips are tangent. 
The number $J$
will remain the same after this assumption. 

We use induction on $S.$ The cases $S=3,4$ we verify directly.
(See Fig. 16 for the maximal possible values of $J$ for $S=3,4$.)
Suppose that the lemma has been proved for
diagrams with fewer than  $S$ strips.
Suppose that a diagram $\Theta$ has $S$ strips. 
Either there is a
strip
that  splits the diagram into two parts 
${\Theta}_1$ and ${\Theta}_2$ (see Fig. 17) with at least two
strips in ${\Theta}_1$ and two strips in ${\Theta}_2,$ or
there is no such strip. If there is no such strip we just draw one more strip
such that is splits the diagram into two parts 
${\Theta}_1$ and ${\Theta}_2$ (see Fig. 17)
and each of these 
parts contains at least two strips. 
Suppose that ${\Theta}_1$ contains $S_1$ strips 
and ${\Theta}_2$
contains $S_2$ strips; in the first case $S=S_1+S_2+1,$ 
in the second case $S=S_1+S_2.$ Apply the induction hypothesis
to the union of ${\Theta}_1$ and the dividing strip. We have
$J_1={\sum}_{k=2}^{S_1+1}N_kk \leq 15(S_1+1-2).$
If we apply induction to the union of ${\Theta}_2$ and the fixed
strip, then we have $J_2={\sum}_{k=2}^{S_2+1}N_kk
\leq 15(S_2+1-2).$ Upon summation,  
$J={\sum}_{k=2}^{S}N_kk\leq J_1+J_2\leq 15 (S_1+S_2-2)\leq 
15 (S-2).$  $\Box$
 
Our goal now is to assign to each $t$-strip  a set of nonsplittable pieces
and $\cal V$-pieces 
in such a way that every piece occurs at most once in the union of these sets.
 
\begin{center}{\bf Dual Forest}\end{center}

 We construct a dual forest of $\Delta$ 
in the following way.  
We plot a vertex of the dual forest in each piece (recall that ``piece''  
means
a nonsplittable piece or a  $\cal V$-piece).

Before defining the of the dual forest of $\Delta$, we consider 
the family $\cal F$ of all subdiagrams of $\Delta$ with the following
property: if $\Theta \in \cal F$ contains some cell of a $t$-strip,
then it contains the whole $t$-strip, and if $\Theta$ contains a cell
of a maximal $G$-subdiagram, then it contains the whole subdiagram. 
We shall define the edges of the  forest for subdiagrams in $\cal F$ 
by induction on the number of $t$-strips.

\begin{df} The {\em piece part of a $t$-strip}, corresponding to a piece $I$,
is the maximal connected set of cells $\tau (I)$ in the $t$-strip 
that border the 
piece $I$ (i.e. having some path in common with the contour of the piece).
\end{df}

\begin{df} Two pieces $I$ and $J$ are called {\em neighbours} if  
$\tau (I)\cap\tau (J)\not=\emptyset$; in other words, if there is a 
cell in a $t$-strip such that its ${\cal U}$-side belongs 
to one piece
and its ${\cal V}$-side to the other.
\end{df}

Suppose the subdiagram $\Theta$ has one $t$-strip (see Fig. 18). 

 We follow one of the sides of the $t$-strip (the $\cal U$-side, for instance)
starting at one of the $t$-edges.
As we meet the first piece $I_1$ (see Fig. 18), 
we draw directed edges from $I_1$ 
 to all the neighbouring
pieces, if there are any (in Fig. 18 these pieces are 
$I_8$ and $I_7$), 
and to each edge we associate the piece part of the $t$-strip
determined by the piece representing the endpoint of the edge.
We also color these associated piece parts.

$S_8$ is assigned to the edge $(I_1, I_8)$, $S_7$ is assigned to the edge
$(I_1, I_7).$ 
(On the Fig. 18 we color $S_7$ and $S_8.$) Then we take the 
next piece $I$ along the $\cal U$-side of the $t$-strip. 
There are three possibilities:
1) $I$ does not have neighbouring pieces;
2) the piece part $\tau (I)$ is already coloured;
3) 
 the piece part $\tau (I)$ is not completely coloured and $I$
has neighbouring pieces.

In the first two cases we go to the next piece, in the third case
we draw directed edges to the vertices in the neighbouring
pieces. And again
to each edge we associate the piece part of the $t$-strip
determined by the piece representing the endpoint of the edge.
We also color these associated piece parts.

We continue this process until we have exhausted all the pieces on the 
 $\cal U$-side of the $t$-strip. 
Finally all the piece parts of the $t$-strip,  associated
to the pieces representing the endpoints of the edges, are colored. 
Every uncoloured cell of the $t$-strip has at least one side on the
boundary of $\Theta .$

 Now suppose that for subdiagrams from $\cal F$ that contain not more
than $s$ strips we have an algorithm to 
construct the Dual forest  and simultaneously colour the piece parts of 
$t$-strips 
assigned 
to the edges 
of the dual forest. Assume furthermore that the algorithm is such that
 the following conditions are satisfied: \begin{enumerate}
\item
The piece parts of the $t$-strip, associated
to the pieces representing the endpoints of the edges, are colored. 
Every uncoloured cell of the $t$-strip has at least one side on the
boundary of the subdiagram.
\item A vertex of the forest cannot be an endpoint of two edges.
\end{enumerate}
Now suppose we have $s+1$ $t$-strips in the subdiagram $\Theta$.
Fix one $t$-strip $S$, such that there are no $t$-strips on one 
side of it. By a simple induction argument such a $t$-strip
always exists. Suppose, for definiteness, that 
there are no $t$-strips on the $\cal U$-side of it. 
Consider the subdiagram ${\Theta}_1$ on the other side ($\cal V$-side)
of this $t$-strip. 
This diagram contains $s$ $t$-strips (see Fig. 19)
We suppose that the Dual forest and the coloring for ${\Theta}_1$ have 
already been constructed and satisfy the induction hypothesis. 
Our purpose is to extend the dual forest to $\Theta$ and to 
colour the piece parts of the $t$-strip $S$.

Consider now the subdiagram ${\Theta}_2$ from $\cal F$ consisting
of $S$ and the neighbouring maximal $G$-subdiagrams. 
We go from the right to the left along the $\cal V$-side of $S,$ 
and repeat the procedure of drawing edges and colouring the 
piece parts of $S$ as in the first step of the induction.

The constructed graph is still a forest, because the graph in 
${\Theta}_1$ is a forest by our induction assumption, and our construction is 
such 
that the new arrows do not produce 
cycles. Moreover, 
all the piece parts of $t$-strips in ${\Theta}$ are assigned
to edges of the dual forest, the uncoloured parts have one side on 
the boundary
of ${\Theta}$ and a vertex of the forest cannot be an endpoint of two edges.
\begin{lm}
\label{main}
Suppose there is  a nonsplittable $k$-piece in  $\Delta$
with the contour label
$$W\equiv W_1X_1W_2X_2\ldots W_kX_k.$$
Let $X_{i_1},\ldots ,X_{i_p}$ be the labels of the  
parts of the boundary of this piece that are on the boundary of
$\Delta.$  Let $n_i$ be the number of $t$-cells corresponding to
the word $W_i.$
Let $M_1, M_2$ be the constants from Corollary 6. Then for any 
$M_3\geq M_2+2c,$ $M_3\geq H,$ 
and for any $i,$ $$\left | n_i\right | \leq
M_1\sum_j \left | X_{i_j}\right | +M_3k.$$
\end{lm}

{\em Proof:} This follows from Corollary \ref{est} and the fact that for a 
nonsplittable piece
the $X_{i_j}$-s that do not belong to the boundary of $\Delta$
are rather short, shorter then $2c,$
hence their sum is less then $2ck.$ $\Box$

\begin{lm} Let $E_1$ be the sum of all cells in all $t$-strips of $\Delta, $
$D$ the length of the boundary of $\Delta, $ and $S$ the number of
$t$-strips (hence $D_1=D-2S$ is the length of the part of the boundary
excluding the $t$-edges). Let $M_4=(1+c\lambda )(15M_3+M_1).$  
Then $E_1\leq M_4D.$ $\Box$
\end{lm}
{\em Proof:} Every $t$-strip is partitioned into piece parts and
parts intersecting the boundary of $\Delta .$ To each piece part we
assign some nonsplittable piece or some maximal $\cal V$-
piece, namely the endpoint of the corresponding edge in the dual forest.
To the remaining parts we assign the intersection with the boundary of $\Delta.$

The label of each part of a $t$-strip that is assigned to a $\cal V$-piece is 
shorter
in $\cal G$ than $2H+R,$ where $R$ is the  
length of the other part of 
the $t$-strip
on the boundary of this $\cal V$-piece, and this part of the $t$-strip
is estimated already not in the $\cal V$-piece but in
nonsplittable pieces (two $\cal V$-pieces cannot be neighbours). 
Recall that $H\leq M_3.$

So, from Lemma \ref{main} we have that 
$$E_1\leq (1+c\lambda )M_1D_1+(1+c\lambda )M_3{\sum}_{k=2}^S N_kk.$$

By Lemma \ref{NK}, $J=\sum_{k=2}^SN_kk\leq 15(S-2).$ 
So $E_1\leq (1+c\lambda )(M_1
D_1+15M_3)=M_4D.$ $\Box$

We apply now

{\bf Transformation 6} Replace each maximal $\cal G$-subdiagram $\Theta$  
in $\Delta$
by a diagram with the same contour, but where the number of cells 
contained in it is not more
than $K$ times the length of the contour of $\Theta$ (it is possible to do 
this, because $\cal G$ satisfies a linear isoperimetric inequality with
the constant $K$).  
\begin{lm}
Let $E$ be the number of cells in  $\Delta ,$
and $M_5=K(1+2M_4c)+M_4.$
Then $$E\leq M_5D.$$
\end{lm}

{\em Proof:} 
If $E_2$ is the total number of cells in all maximal ${\cal G}$-subdiagrams,
then
$$E_2\leq K(D_1+2M_4cD).$$
Finally 
$$E=E_1+E_2\leq M_4D+K(D_1+2M_4cD)=M_5D.\Box$$

The Proposition and Theorem 1 are proved.$\Box $

{\em Proof of Corollary 1}

Let ${\cal G}$ be a hyperbolic group,
$\cal A$ and $\cal B$ isomorphic virtually cyclic subgroups, hence (\cite{G},
\cite{Ol}) quasiisometrically embedded. Then
the separated HNN-extension
${\cal H}=<{\cal G},t\mid  {\cal A}^t={\cal B}>$
is hyperbolic by Theorem 1.

Now suppose that the HNN-extension $\cal H$  is not separated.

First, suppose that there is an element $g\in \cal G$ such that
${\cal A}^g\cap {\cal B}$ is infinite. Then there is an element
$b=a^g$ of infinite order, where $b\in \cal B$ and $a\in\cal A.$
Then, for some nonzero integers $m$ and $n,$ we have $t^{-1}a^nt=b^{m},$
hence $ t^{-1}a^nt=g^{-1}a^{m}g$ and
$$(tg^{-1})^{-1}a^ntg^{-1}=a^{m}.$$

It follows from \cite{G} (Corollary 8.2.c) that if, in a hyperbolic group,
an element $y$ has infinite order, and $m$ and $n$ are nonzero integers,
then the equation $x^{-1}y^nx=y^{m}$ implies that the subgroup
generated by $x$ and $y$ is virtually cyclic. In particular,
$|m|=|n|$ and we say that $x$ nearly commutes with $y^m.$

But in our case both elements $tg^{-1}$ and $a$ have infinite order and
they are not powers of the same element (because the reduced form
of $(tg^{-1})^{k_1}$ in the HNN-extension $\cal H$ is different from the
reduced form $a^{k_2}$ for any $k_1, k_2$).
Hence the subgroup generated by $tg^{-1}$ and $a$ cannot be
virtually cyclic.

Suppose now that neither $\cal A$ nor $\cal B$ is conjugate separated.
Then, for some $g_1\in {\cal G}\backslash {\cal A}$ and
$g_2\in {\cal G}\backslash {\cal B},$ both sets
$$S_1=\{a\in {\cal A} \mid g_1^{-1}ag_1\in {\cal A}\}$$
and 
$$S_2=\{b\in {\cal B} \mid g_2^{-1}bg_2\in {\cal B}\}$$
are infinite. Since $S_1$ and $S_2$ are then infinite subgroups of the
virtually cyclic groups ${\cal A}$ and $\cal B,$ there are  elements $c\in S_1$ 
and $d\in S_2$ of infinite
order. The inclusion $g_1^{-1}
cg_1\in {\cal A}$ implies that $g_1$ nearly commute with a power of $c.$
Also, $g_2^{-1}dg_2\in {\cal B}$ implies that $g_2$ nearly commutes with
a power of $d.$ If $\cal H$ were hyperbolic 
we would be able to find numbers $m$ and $n$ such that
$g_1^{-2}c^mg_1^2=c^m,$ $g_2^{-2}d^ng_2^2=d^n$ 
and $t^{-1}c^mt=d^{n}.$ It is easy to see that
the subgroup $<(tg_2^2t^{-1}g_1^2)^{2}, c^m>$ is free 
abelian, a contradiction. $\Box$
\section{Free products with amalgamation}

The proof of Theorem 2 is quite similar and not as complicated as the
proof of Theorem~1, so this section will be quite brief, 
more in the vein of a guided exercise than 
a proof.

Let ${\cal G}_1$, ${\cal G}_2$ and $\cal U$, $\cal V$ be 
as in Theorem 2, ${\cal R}={\cal G}_1\ast_{\cal U=\cal V}{\cal G}_2.$

From the hyperbolicity  of ${\cal G}_1$ and ${\cal G}_2$ 
it follows that there is a 
constant $K$ such that for any element $W_1=1$ (resp. $W_2=1$) 
in ${\cal G}_1$ (resp. ${\cal G}_2$) there
is a reduced diagram over ${\cal G}_1$ (resp. ${\cal G}_2$) 
with boundary label $W_1$ (resp. $W_2$) and number of cells
$\leq K\left |\left | W_1\right |\right | $ ($\leq K\left |\left | W_2\right 
|\right | $). 

Let $A=\{a_1,\ldots ,a_n\}$ and $B=\{b_1,\ldots ,b_n\}$ 
be the distinguished generating sets 
for ${\cal G}_1$ and ${\cal G}_2$ respectively
, such that $\phi (a_i)=b_i.$
The contents of this section is the proof of the following 
\begin{prop} There is a linear function $L_1$ of a single variable, 
depending only on 
${\cal G}_1,{\cal G}_2,
A,B,$ such that
for any element $W=1$ in $\cal R$ there
is a diagram over $\cal R$ with boundary label $W$ and number of cells
$\leq L_1(\left |\left | W\right |\right | ).$
\end{prop}

A $(\cal U,\cal V)$-cell is a cell with contour $a_ib_i^{-1}.$
In this section $(\cal U,\cal V)$-cells will play  the role of $t$-cells.

Let $\Delta$ be a minimal diagram over $\cal R$ 
with boundary label $W$ (this means that $\Delta$ contains a minimal
possible number of $(\cal U,\cal V)$-cells).

A $(\cal U,\cal V)$-strip is a subdiagram, with boundary 
$a_{i_1}\ldots a_{i_k}b_{i_k}^{-1}\ldots b_{i_1}^{-1}$ consisting
of $(\cal U,\cal V)$-cells, that begins and ends on the boundary of
the diagram $\Delta,$ and is minimal with this property. 
$(\cal U,\cal V)$-strips will play the role of 
$t$-strips. A $(\cal U,\cal V)$-cell and $(\cal U,\cal V)$-strip 
are shown in Fig. 20.

\begin{lm} 
A diagram $\Delta$ cannot contain a $(\cal U,\cal V)$-annulus.
\end{lm}
The proof follows from the minimality of $\Delta$. $\Box$

The diagram $\Delta$ consists of maximal 
${\cal G}_1$- and ${\cal G}_2$-subdiagrams 
that are glued to each 
other through $(\cal U,\cal V)$-strips.
A typical form of $\Delta$  is shown on Fig. 21; the same diagram
is schematically shown in Fig. 22 (the notion of
a $(\cal U,\cal V)$-strip is clear from Fig. 22). The
partition of the set of all $(\cal U,\cal V)$-cells into
$(\cal U,\cal V)$-strips is not necessary  unique; we just take
some partition.

The notions of an island, 
nonsplittable piece (in a maximal ${\cal G}_1$-
or ${\cal G}_2$-subdiagram) and of $\cal V$-piece are the same as 
in the previous section.

We perform Transformations 2 and 3 
on $\Delta$. 
 After each transformation we reduce ${\Delta}$ to a minimal diagram.
 
\begin{lm}
\label{mng}
Every maximal ${\cal G}_1$-subdiagram can be transformed into  
a disjoint
union of nonsplittable pieces 
with boundary labels 
of the form
$$ W_{1}X_1W_{2}X_2\ldots W_{k}X_k=1,$$
where all the $X_i$ correspond to boundary paths.
\end{lm}
The proof is a simpler version of the proof of Lemma \ref{na},
so we omit it.

We construct a dual forest as we did in the previous section,
but instead of $\cal V$-pieces and 
nonsplittable pieces in ${\cal G}_2$-subdiagrams
 we just use maximal ${\cal G}_2$-
subdiagrams.

\begin{lm}
\label{main1}
Suppose that  ${\Delta}$ contains a  
nonsplittable ${\cal G}_1$-piece
with the contour label
$$W_{1}X_1W_{2}X_2\ldots W_{k}X_k.$$

Let $n_i$ be the number of $(\cal U, \cal V)$-cells corresponding to $W_i.$
Then for $M_6=M_1+M_2$ (these are the constants from Corollary 6)  
and for any $i,$ $$\left | n_i\right | \leq
M_6\sum_i \left | X_{i}\right | .$$
\end{lm}

Proof. All the $X_i$-s are nonempty words (recall that if $X_i$=1, then
we consider the union of two strippaths $W_i$ and $W_{i+1}$ as one
strippath). Hence $\sum_i \left | X_{i}\right |\geq k.$ $\Box$
\begin{lm} Let $E_1$ be the sum of all the cells in 
all the $(\cal U,\cal V)$-strips of 
${\Delta}, $ $M_7=M_6+c\lambda M_6,$
and let $D$ be the length of the boundary of ${\Delta}.$ 
Then $E_1\leq M_7D.$
\end{lm}
{\em Proof:} Every $(\cal U,\cal V)$-strip is subdivided into parts such that 
each 
part is 
either assigned to some nonsplittable ${\cal G}_1$-piece, or 
to a path on the boundary of ${\Delta}$
not included in the boundary of some nonsplittable 
${\cal G}_1$-piece, or to a maximal ${\cal G}_2$-subdiagram. 
Each nonsplittable ${\cal G}_1$-piece and each such path on 
the boundary of ${\Delta}$
cannot be assigned to more then one $(\cal U,\cal V)$- strip and to more  
then one distinct part of  
this strip.
The length of each piece of a $(\cal U,\cal V)$-strip, assigned 
a maximal ${\cal G}_2$-subdiagram $\Theta,$ is not greater in $\cal G$ 
than the sum of lengths of all other
parts of $(\cal U,\cal V)$-strips in $\Theta$ (that are assigned
to some
nonsplittable  ${\cal G}_1$-pieces or to parts of the boundary of 
$\Delta$) 
plus the length of all the boundary-paths of $\Theta.$

Suppose $M_6\geq 1.$ From Lemma \ref{main1} we have that 
$$E_1\leq M_6D+c\lambda M_6D=M_7D.\Box$$

Now we apply to $\Delta$ the analog of Transformation 6. We replace
all maximal ${\cal G}_i$-subdiagrams by diagrams with the same contour,
but where the number of cells
is less then the length
of the contour times $K.$
\begin{lm}\label{18}
Let $E$ be the number of cells in  ${\Delta} ,$
and $M_8=K(1+2M_7c)+M_7.$
Then $$E\leq M_8D.$$
\end{lm}

{\em Proof:} 
If $E_2$ is the sum of the cells in all maximal ${\cal G}_i$-subdiagrams
then
$$E_2\leq K(1+ 2M_7c)D.$$
Finally, 
$$E=E_1+E_2\leq M_7D+K(1+2M_7cD)=M_8B. \ \ \Box$$
The Proposition and Theorem 2 now follow from  Lemma \ref{18}.
 
{\em Proof of Corollary 2}

In one direction, Corollary follows from Theorem 2 and the fact that
virtually cyclic subgroups of hyperbolic groups are quasiisometrically
embedded.

Suppose now that $\cal A$ is not conjugate separated in ${\cal G}_1$ 
and $\cal B$ is not conjugate separated in ${\cal G}_2.$
Then for some $g_1\in {\cal G}_1\backslash {\cal A}$ and
$g_2\in {\cal G}_2\backslash {\cal B}$ both sets
$$S_1=\{a\in {\cal A} \mid g_1^{-1}ag_1\in {\cal A}\}$$
and 
$$S_2=\{b\in {\cal B} \mid g_2^{-1}bg_2\in {\cal B}\}$$
are infinite. Since $S_1$ and $S_2$ are then infinite subgroups of the
virtually cyclic groups ${\cal A}$ and $\cal B,$ there are  elements of infinite
order $c\in S_1$ and $d\in S_2.$ The inclusion $g_1^{-1}
cg_1\in {\cal A}$ implies that $g_1$ nearly commutes with a power of $c.$
Also, $g_2^{-1}dg_2\in {\cal B}$ implies that $g_2$ nearly commutes with
a power of $d.$ There is a common power $z$ of $c$ and $d$ 
such that
$g_1$ and $g_2$ both nearly commute with $z.$ It is easy to see that
the subgroup $<(g_1g_2)^{2}, z>$ is free abelian, a contradiction. 
$\Box$
\section{Other sufficient conditions}

Let $\cal H$ be the fundamental group of a finite graph $\Gamma $ of
groups (relative to some maximal subtree $T$ of $\Gamma$ ) with vertex groups
${\cal G}_i, i=1,\ldots p,$ 
edge groups ${\cal U}_{ij},$ such that
${\cal U}_{ij}\leq {\cal G}_i,$ and embeddings 
$\tau:{\cal U}_{ij}\rightarrow  
{\cal G}_j$, such that $\tau({\cal U}_{ij})={\cal V}_{ij}={\cal U}_{ji}\leq 
{\cal 
G}_j.$
Then $\cal H$ is generated by the groups ${\cal G}_i$
and additional elements $t_{ij},$ which 
are in bijective correspondence with
the non-$T$ edges. $\cal H$ has, in addition to the relations of groups
${\cal G}_i,$ the following defining relations:
$u=\tau (u)$ for any $u\in {\cal U}_{ij},$ with $(ij)$ is a $T$-edge,
and $u^{t_{ij}}=\tau (u)$ for all $u\in{\cal U}_{ij}$ with $(ij)$
a non $T$-edge.

\begin{df} A diagram of the type shown in Fig. 23
is called an {\em $h$-rectangular} subdiagram over $\cal H$, if the following
conditions are satisfied:
\begin{enumerate}
\item The strips are either 
$({\cal U}_{ij},{\cal V}_{ij})$-strips or 
$t_{ij}$-strips (which we will also call 
$({\cal U}_{ij},{\cal V}_{ij})$-strips).
\item The subdiagrams between strips are maximal ${\cal G}_i$-subdiagrams.
\item The boundary-paths $p_k, q_k$ are shorter than some fixed number $h.$
\item If the labels of the strippaths of a maximal ${\cal G}_i$-subdiagram
belong to the same edge group ${\cal U}_{ij},$  
and both strips that bound it are
$({\cal U}_{ij},{\cal V}_{ij})$-strips, then the label of at least one of its
two boundary paths does not belong to ${\cal U}_{ij}$.
\end{enumerate}

The number $n$ of $({\cal U}_{ij},{\cal V}_{ij})$-strips in the diagram is 
called
the {\em length} of the diagram; the paths $[x_1,x_2],\  [x_3,x_4]$
are called the {\em sides} of the diagram. The lengths of the labels
of the sides
are taken in the corresponding edge groups. If the labels of the two boundary
paths
of each maximal ${\cal G}_i$-subdiagram are the same 
($\phi (p_k)=\phi (q_k)$ ), 
then the diagram is called {\em a conjugacy $h$-rectangular diagram.}
\end{df}

We obtain sufficient conditions for the hyperbolicity of $\cal H$ 
which are weaker than
the conditions in Theorems 1 and 2.
\begin{th}
Let $\cal H$ be the fundamental group of a finite graph of
groups, with the edge groups ${\cal U}_{ij}$ quasiisometrically embedded in
the corresponding vertex groups ${\cal G}_i$ and ${\cal G}_j$ 
($\varepsilon$-quasiconvex). 
Suppose that all the vertex groups
${\cal G}_i$ are hyperbolic, and $\delta $ is the maximum of hyperbolicity
constants of the vertex groups. 
Let $H=8\delta + \varepsilon .$ 
If there exists a number $n$ such that only a finite number of elements
in $\cal H$ can be labels of the sides of a  reduced 
conjugacy $2H$-rectangular
diagram of length $n$, then $\cal H$ is hyperbolic.\end{th}

The condition of the theorem implies that all the elements in this
finite set have finite order.

We will prove the theorem after formulating the 
following corollaries.

\begin{cy}
Let $\cal H$ be a fundamental group of a finite graph of
groups, with edge groups ${\cal U}_{ij}$ quasiisometrically embedded in
the corresponding vertex groups ${\cal G}_i$ and ${\cal G}_j$ 
Suppose that all the vertex groups
${\cal G}_i$ are torsion-free hyperbolic 
, and $\delta $ is the largest of the constants of hyperbolicity
of the vertex groups. 
Let $H=8\delta + \varepsilon .$ 
If there exists a number $n$ such that there are no reduced 
conjugacy $2H$-rectangular
diagrams of length $n$, then $\cal H$ is hyperbolic.\end{cy}

\begin{cy} Let ${ {\cal G}}_1$, ${ {\cal G}}_2$ be hyperbolic groups, 
${ {\cal U}}\leq { {\cal G}}_1$, ${ {\cal V}}\leq 
{ {\cal G}}_2$ quasiisometrically embedded, and $\phi :{\cal U}\rightarrow {\cal 
V}$
an isomorphism. Suppose that  
there exists a number $n$ such that the set
$$h_n\ldots (g_2(\phi ^{-1}(h_1(\phi(g_1{\cal U}g_1^{-1}\cap {\cal U}))
h_1^{-1}\cap {\cal V}))g_2^{-1}\cap {\cal U}
)\ldots h_n^{-1}
\cap {\cal V}$$ is finite (here all $g_i\in {\cal G}_1\backslash {\cal U}, 
h_i\in 
{\cal G}_2\backslash {\cal V}$.)
Then the group ${ {\cal G}}_1\ast_{{ {\cal U}}={ {\cal V}}}{ {\cal G}}_2$
is hyperbolic.
\end{cy}
For HNN-extensions $<{\cal G},t | {\cal U}^{t}=\phi ({\cal U}) = {\cal V}>$ 
there 
is a more
complicated condition : 
\begin{cy}
Let ${ {\cal G}}$ be a hyperbolic group, 
${ {\cal U}}\leq { {\cal G}}$, ${ {\cal V}}\leq 
{ {\cal G}}$ quasiisometrically embedded, and $\phi :{\cal U}\rightarrow {\cal 
V}$
an isomorphism. Suppose  
there exists a number $n$ such that for any
$C_1,..., C_n \in \{{\cal U},{\cal V}\}$ the set
$$(g_n\ldots (\phi ^{\alpha_3}(g_2(\phi^{\alpha_2}(g_1C_1g_1^{-1}\cap C_2))
g_2^{-1}\cap C_3))\ldots g_n^{-1})\cap C_{n+1}$$   is finite.
Here, if $C_i=C_{i+1}$ then $g_i\not\in C_i;$
$\alpha _i=1$ if $C_i={\cal U}$ and $\alpha _i=-1$ if $C_i={\cal V}$.  

Then the group ${\cal H}=<{\cal G},t | {\cal U}^{t}=\phi ({\cal U}) = 
{\cal V}>$
is hyperbolic.
\end{cy}

{\em Proof of the theorem:}

We will prove that any diagram over $\cal H$ satisfies a linear
isoperimetric inequality. The idea behind the proof is exactly the
same as the idea of the proof of Theorems 1 and 2. 

It can be shown, as in the proof of Lemma 4, that if the sides of any
minimal conjugacy 
$2H$-rectangular diagram of length $n$ are bounded by a constant $C,$
then there is a number $\bar C$ such that the sides of 
any minimal $2H$-rectangular
diagram of length $n$ are bounded by $\bar C.$ (If they are not bounded, then 
one can find infinitely many elements that are the labels of the sides of 
conjugacy $2H$-rectangular diagram.)
Let $\Delta$ be a minimal diagram over $\cal H$. It is subdivided by
$({\cal U}_{\alpha\beta}, {\cal V}_{\alpha\beta})$-strips 
into maximal ${\cal G}_{\alpha}$-subdiagrams. 
The contour of each ${\cal G}_{\alpha}$-subdiagram is a word 
$$W_1X_1W_2X_2\ldots W_kX_k,$$
where $W_1,\ldots ,W_k\in \{{\cal U}_{\alpha\beta}, {\cal V}_{\gamma\alpha}|
\alpha,\beta,\gamma \in\{1,\ldots ,p\}\}.$

We shall give a slightly
different definition of a nonsplittable decomposition, than in Section 3.
First, the following 
\begin{df}
Given a decomposition 
$$W\equiv W_1X_1W_2X_2\ldots W_kX_k,$$
where $W_1,\ldots ,W_k\in \{{\cal U}_{\alpha\beta}, {\cal V}_{\gamma\alpha}|
\alpha,\beta,\gamma =1,\ldots p\},$
$W_i$ are geodesic in the corresponding groups, $X_i\in {\cal G}_{\alpha}$ 
are  reduced, and if $W_i,W_{i+1}\in 
{\cal U}_{\alpha\beta} (\cal V_{\gamma\alpha}), $ 
then $X_j\not =1$,  
we call the decomposition a {\em ${\cal G}_{\alpha}$-edges decomposition}.

If $W=1,$ and the indices in the above definition are taken modulo k, 
then it is called a {\em cyclic ${\cal G}_{\alpha}$-edges decomposition.}
\end{df}
\begin{df} A cyclic ${\cal G}_{\alpha}$-edges decomposition
is called {\em splittable} if one of the following holds
\begin {enumerate}
\item There is a $j\leq k$ such that 
$W_1\equiv W_{11}W_{12},\  W_{11},W_{12}\in { {\cal U}_{\alpha\beta}} ,$
$W_j\equiv W_{j1}W_{j2}, \ W_{j1},W_{j2}\in {{\cal U}_{\alpha\beta}} $
 and 
$$W_{12}X_1W_2X_2\ldots W_{j1}=W'\in { {\cal U}_{\alpha\beta }},$$
where $$\left | W_{12}\right | _{{\cal U}_{\alpha\beta }} +\left | W_{j1}\right 
| 
_{{\cal U}_{\alpha\beta }}>
\left | W'\right | _{{\cal U}_{\alpha\beta }}+\left | W_{j1}W'^{-1}W_{12}
\right |  _{{\cal U}_{\alpha\beta }}.$$   
\item There is a $j\leq k$ such that 
$W_1\equiv W_{11}W_{12},\  W_{11},W_{12}\in { {\cal V}_{\gamma\alpha}} ,$
 $w_j\equiv w_{j1}w_{j2}, \ w_{j1},w_{j2}\in {{\cal V}_{\gamma\alpha}} $
 and 
$$W_{12}X_1W_2X_2\ldots W_{j1}=W'\in { {\cal V}_{\gamma\alpha }},$$
where $$\left | W_{12}\right | _{{\cal V}_{\gamma\alpha }}
 +\left | W_{j1}\right | _{{\cal V}_{\gamma\alpha }}>
\left | W'\right | _{{\cal V}_{\gamma\alpha}}
+\left | W_{j1}W'^{-1}W_{12}
\right |  _{{\cal V}_{\gamma\alpha }}.$$
\item There is a $j\leq k$ such that $W_1\equiv W_{11}W_{12}W_{13}, 
W_{11}\not =1,
W_{13}\not =1 , W_{11},W_{12},W_{13}\in {\cal U}_{\alpha\beta} 
({\cal V}_{\gamma\alpha}),$
 $W_j\equiv W_{j1}W_{j2}W_{j3}, W_{j1}\not =1, W_{j3}\not = 1
W_{j1},W_{j2},W_{j3}\in {\cal U}_{\alpha{\beta}_1} 
({\cal V}_{{\gamma}_1\alpha})$
and $$W_{13}X_1W_2X_2\ldots W_{j1}=T,$$
$$W_{j3}X_j\ldots W_{11}=T_1,$$

where  $\left | T\right |, \left | T_1\right | < H.$ 
In this case we always take pairs $W_{11},W_{j3}$
and $W_{13}, W_{j1}$ of minimal length (in the corresponding edge groups) 
among the pairs with the same property. (The length of the pair
 $W_{k},W_{t}$ is not less than the length of the pair
$\bar W_{k}, \bar W_{t}$ if the length of $ W_k$ is not less than
the length of $\bar W_k $
and the length of $ W_t$ is not less than
the length of $\bar W_t. $
)

\end {enumerate}

A cyclic decomposition of $W$ is called 
{\em cyclically nonsplittable} if all   
the permutations
of the form 
$$W_iX_i\ldots W_kX_kW_{k+1}\ldots W_{i-1}X_{i-1}=1.$$
are nonsplittable.
\end{df}

As in Section 3, a maximal nonsplittable ${\cal G}_i$-subdiagram 
is called a {\em nonsplittable piece}.
\begin{df} A {\em thin bridge} is a ${\cal G}_i$-subdiagram having 
boundary label 
$$W_1X_1W_2X_2$$
where the $X_i$'s are shorter than $H$ (the $X_i$'s can be trivial) 
and the $W_i$'s are 
labels of strip-paths and geodesic in the corresponding edge groups; 
moreover, the thin bridge is required to be maximal
among such ${\cal G}_i$-subdiagrams.
\end{df}
\begin{lm} 
Every maximal $\cal \cal G$-subdiagram in $\Delta$ 
consists of nonsplittable pieces connected 
by thin bridges.
\end{lm}
The proof is very similar to the proof of Lemma 9.
The diagram $\Delta$ is minimal hence cannot contain 
a $({\cal U}_{ij}, {\cal V}_{ij})$-
annulus.
One can perform on $\Delta$ the obvious analogs of 
Transformations 2 and 3.

We can now construct a dual forest in $\Delta$ 
like we did for $HNN$-extensions,
but instead of nonsplittable pieces 
and ${\cal V}$-pieces we use
nonsplittable pieces 
and thin bridges.

\begin{lm}
\label{est1} There are constants $M_1$, $M_2$ such that if,
in a vertex group ${{\cal G}_{\alpha}},$ 
we have a cyclically nonsplittable 
decomposition of the word 
$$W\equiv W_1X_1W_2X_2\ldots W_kX_k=1,$$
then, for
 $W_i\in {\cal U}_{\alpha \beta} ({\cal V}_{\gamma\alpha}),$ we have
${\left | W_i\right | }_{{\cal U}_{\alpha \beta} ({\cal V}_{\gamma\alpha})}
\leq M_1\sum_i \left | X_{i}\right | + M_2k$.
\end{lm}
The proof is similar to the proof of Corollary \ref{est}.

The proof of the following lemma repeats the proof of Lemma \ref{NK}.
\begin{lm} 
\label{NK1}
Let $N_k$ be the number of nonsplittable $k$-pieces if $k\geq 3,$ and 
let $N_2$
be the number of nonsplittable $2$-pieces plus the number of thin bridges
in
a diagram $\Delta $ over the group ${\cal H}.$ Let $S$ be the number of
$({\cal U}_{ij}, {\cal V}_{ij})$-strips in it, $S\geq 3.$ 
Then  $J=\sum_{k=2}^SN_kk\leq 15 (S-2).$
\end{lm}

Notice that $S$ is no larger than $p$ (the number of vertex groups) times
the sum of the number of distinct boundary-paths of $\Delta$ and 
the number of $t_{ij}$-edges
on the boundary of $\Delta $; hence, if $D$ is the
length of the boundary of $\Delta,$ then $S\leq pD.$

The proof of the following lemma repeats the proof of Lemma \ref{main}.

\begin{lm}
\label{main2}
Suppose we have a nonsplittable $k$-piece  in $\Delta$
with the contour label
$$W\equiv W_1X_1W_2X_2\ldots W_kX_k.$$
 
Let $X_{i_1},\ldots ,X_{i_p}$ be pieces of the boundary of 
$\Delta.$  Let $n_i$ be the number of 
$({\cal U}_{\alpha \beta},{\cal V}_{\alpha\beta})$-
cells corresponding to
the word $W_i.$
Then, for  any $M_3\geq M_2+2c$ 
and for any $i,$ we have $$\left | n_i\right | \leq
M_1\sum_j \left | X_{i_j}\right | +M_3k.$$
\end{lm}

In constructing the dual forest we assigned to each piece part of each 
$({\cal U}_{\alpha \beta},{\cal V}_{\alpha\beta })$-strip 
an edge of the dual forest. This edge is associated  
either with some nonsplittable piece
or with some thin bridge, containing the endpoint of the edge. 

The length of those parts of the 
$({\cal U}_{\alpha \beta},{\cal V}_{\alpha\beta })$-strips
that are assigned to the edges having endpoints in nonsplittable pieces,
are estimated in these pieces. The only problem is to estimate
the length of the 
$({\cal U}_{\alpha \beta},{\cal V}_{\alpha\beta })$-strips
assigned to edges having endpoints in thin bridges.

Consider the configuration of thin bridges as shown in Fig 24.
The arrows correspond to the edges of the dual forest. 

The sides of 
any $2H$-rectangular
diagram of length $n$ are bounded by $\bar C$.
Take the part of the  
$({\cal U}_{\alpha \beta},{\cal V}_{\alpha\beta })$-strip $S$
assigned to the piece associated with the endpoint of the first edge
($I_1$). $S$ can be split into parts, such that for each
part one of the following possibilities applies : \begin{enumerate} 
\item we can
construct a $2H$-rectangular subdiagram of length more than $n$ 
starting with this part; \item 
the length of the rectangular subdiagram connecting 
this part of the strip $S$ with some part of a strip 
$R,$ assigned to an edge of the dual forest with 
endpoint in some nonsplittable piece, is less than $n$ (then 
the length of this part can be bounded by the length of
$R$ times some constant $C_1$, depending only on $H$ and $n$).
\end{enumerate}

Now, taking $M_3>\bar C$ and $M_4=nC_1M_1+15npM_3,$ we have that the number of
all  $({\cal U}_{\alpha \beta},{\cal V}_{\alpha\beta })$-cells in $\Delta$
is less than $M_4D,$ where $D$ is the length of the boundary of
$\Delta .$ And, taking $M_5$ as in Lemma 13, we get the
linear isoperimetric inequality for $\Delta$ with the constant $M_5.$
The Theorem has been proved.$\Box$

To prove the Corollaries it suffices to notice that the 
conditions in the statements of the Corollaries
imply the conditions of the theorem. 

\section{Some results on quasiconvexity}

A subset $Y$ in a geodesic space $\Gamma$ is called quasiconvex for some
$\epsilon \geq 0$ if every geodesic segment $[y_1,y_2]$ with  endpoints in 
$Y$ lies $\epsilon$-close to $Y.$ A subgroup $U$ of a group $G$ is called 
{\em quasiconvex} if ${\Gamma}(U)$ is quasiconvex in ${\Gamma}(G).$

In this section we prove the following theorems.
\begin{th}
\label{qc}
Let ${\cal H}=<{\cal G},t |{\cal U}^t={\cal V}>$ be hyperbolic
with $\cal U$ 
quasiconvex in ${\cal H}.$
Then ${\cal G}$ is quasiconvex in ${\cal H},$ and hence hyperbolic.
\end{th}

\begin{th}
\label{qc1}
Let ${\cal H}$ be a separated HNN-extension
${\cal H}=<{\cal G},t |{\cal U}^t={\cal V}>$ 
with $\cal G$ hyperbolic, 
${\cal U}$ and 
${\cal V}$ quasiconvex in ${\cal G}.$ Then ${\cal G}$
is quasiconvex in ${\cal H}.$
\end{th}
{\em Proof of Theorem \ref{qc}:} 
We have to show that there exist a number $\lambda $
such that the length of arbitrary geodesic in $\cal G$ is 
shorter than $\lambda$ times the length of the same element in $\cal H.$
Let $W$ be a geodesic word in $\cal G$
and $V$ a geodesic word in $\cal H$ such that $V=W$ in $\cal H.$
Let $p$ be a path such that $\phi (p) =W$ and $q$ be a path such that
$\phi (q)=V.$ The subgroup $\cal U$ is quasiconvex in $\cal H,$ hence 
$\cal V$ is quasiconvex in $H$. 
Let $\lambda$ be a number such that every geodesic in $\cal U$ or$\cal
V$ is $\lambda$-quasigeodesic in $\cal H.$ 
Let $\Delta $ be a minimal diagram over $\cal H$ 
with the boundary $qp^{-1}.$ 
Then a typical form of $\Delta$ is shown in Fig. 25. The path $p$
is shorter than the path $s$ in Fig. 25, but the path $s$ is shorter
than $\lambda |q|.$ 

The theorem is proved.

Exactly the same reasoning can be used to prove the following more general
result.

Suppose we have a finite graph of groups, 
with finitely generated   edge groups, and the fundamental group ${\cal H}$
   of the graph is hyperbolic. It then follows that if for some vertex group
    ${\cal G}(v)$ all incoming edge groups are quasi-convex in the whole  
    group ${\cal H}$ , then the vertex group ${\cal G}(v)$ 
is quasi-convex in ${\cal H}.$  The result 
in this formulation was obtained
by I. Kapovich (who used a different technique) \cite{IK}.

{\em Proof of Theorem \ref{qc1}:}
We will show that if $L$ is a linear function as in
Proposition~2, then the length of an arbitrary geodesic in $\cal G$ is 
shorter than $L(\ell),$ where $\ell$ is the length of the same element 
in $\cal H.$
Let $W$ be a geodesic word in $\cal G$
and $V$ a geodesic word in $\cal H$ such that $V=W$ in $\cal H.$
Let $p$ be a path such that $\phi (p) =W$ and $q$ a path such that
$\phi (q)=V.$ Let $\Delta $ be a minimal diagram over $\cal H$ 
with the boundary $qp^{-1}.$ 
Then a typical form of $\Delta$ is shown in Fig. 25. We construct
the dual forest for $\Delta$ starting from the pieces between
$p$ and $s$ as shown in Fig. 26. Then $s$ is shorter than $L(|q|)$
and $p$ is shorter than $s.$ $\Box$

The following theorems can be proved by a similar technique.
\begin{th}
\label{qc2}
Let ${\cal H}={\cal G}_1\ast_{\cal U}{\cal G}_2$ be 
a hyperbolic group,  with ${\cal U}$  
quasiconvex in ${\cal H}.$ Then ${\cal G}_1$ and
${\cal G}_2$
are quasiconvex in ${\cal H},$ and hence hyperbolic.
\end{th}
\begin{th}
\label{qc3}
Let  
${\cal H}={\cal G}_1\ast_{{\cal U}={\cal V}}{\cal G}_2,$
with ${\cal G}_1$ and ${\cal G}_2$ hyperbolic,
${\cal U}$ quasiconvex in ${\cal G}_1,$  
${\cal V}$ quasiconvex in ${\cal G}_2$ and $\cal U$ conjugate separated
in ${\cal G}_1.$ 
Then ${\cal G}_1$ and ${\cal G}_2$
are quasiconvex in ${\cal H}.$ 
\end{th}
In the situation where $\cal U$ is malnormal in ${\cal G}_1$ and 
$\cal V$ is malnormal in ${\cal G}_2,$ this result can also be deduced from
\cite{PP}.
 \section{Applications to exponential groups}

Let $A$ be an arbitrary associative ring with identity and $G$ a group.
Fix an action of the ring $A$ on $G$, i.e. a map $G \times A \rightarrow
G$. The result of the action of $\alpha \in A$ on $g \in G$ is written
as $g^\alpha$. Consider the following axioms:

\begin{enumerate}
\item $g^1=g$, $g^0=1$,  $1^\alpha = 1$ ;
\item $g^{\alpha +\beta}=g^\alpha \cdot g^\beta, \ g^{\alpha \beta} =
(g^\alpha)^\beta$;
\item $(h^{-1}gh)^\alpha = h^{-1}g^\alpha h;$
\item $[g,h]=1 \Longrightarrow (gh)^\alpha = g^\alpha h^\alpha.$
\end{enumerate}

\begin{df}
Groups with $A$-actions satisfying axioms 1)--4) are called {\em $A$--groups}. 
\end{df}
In particular, an arbitrary group $G$ is a ${\bf Z}$-group. We now recall
the definition of $A$-completion in the case where $A$ is the field 
of rationals ${\bf Q}$. Such completions were studied by G. Baumslag in
\cite{bau}, \cite{baum1}. $A$-completions for arbitrary rings $A$
were investigated in \cite{MyasExpo2}. We will use some results 
and constructions from the latter article. 
\begin{df}
Let $G$ be a group . 
Then a ${\bf Q}$--group $G^{\bf Q}$ 
together with a homomorphism 
$G\rightarrow G^{\bf Q}$
is called  a tensor ${\bf Q}$--completion
of the group $G$ if $G^{\bf Q}$ satisfies the following universal
property:
for any ${\bf Q}$--group $H$ and a homomorphism 
$\varphi: G \rightarrow H$
there exists a unique ${\bf Q}$--homomorphism $\psi: G^{\bf Q} \rightarrow H$
(a homomorphism that commutes with the action of ${\bf Q}$) 
such that the following diagram commutes:

\medskip

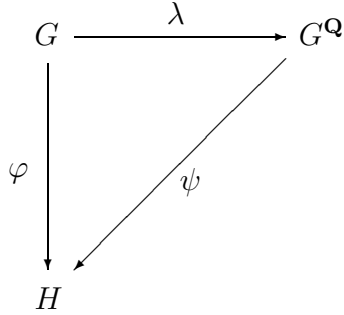
\begin{center}

\begin{picture}(100,100)(0,0)
\put(0,100){$G$}
\put(100,100){$G^{\bf Q}$}
\put(0,0){$H$}
\put(15,103){\vector(1,0){80}}
\put(5,93){\vector(0,-1){78}}
\put(95,95){\vector(-1,-1){80}}
\put(-10,50){$\varphi$}
\put(55,45){$\psi$}
\put(50,108){$\lambda$}
\end{picture}

\end{center}

\medskip
\end{df}
It was proved in \cite{bau} (see also \cite{MyasExpo2}) 
that for every group $G$  the tensor ${\bf Q}$-completion of $G$ 
exists and is unique. 

   In this section we describe the ${\bf Q}$-completion $G^{\bf Q}$ 
of a torsion-free
hyperbolic group $G$ as the union of an effective chain of hyperbolic subgroups.
This allows one to apply techniques of hyperbolic group
theory  to solve various algorithmic
problems in $G^{\bf Q}$, in particular, to construct effectively 
some natural  normal forms for its elements (induced
by the normal forms of amalgamated 
free products). 
 
First of all, let us describe the construction of 
the {\it complete tensor 
extension of centralizers of an arbitrary torsion-free hyperbolic group 
$G$ by the ring ${\bf Q}$} (see  \cite{MyasExpo2}).

Let $C= C_G(v)=C(v)$ be a centralizer in $G$ and $v$ not a proper power, 
i.e. $C(v)=<v>$. The $\bf Q$-extension of the 
centralizer $C$ is by definition a free product with amalgamation
$$G(C, {\bf Q})= G\ast_C{\bf Q},$$
where $C\simeq\bf Z\leq \bf Q.$ The group $G(C,\bf Q)$ can be obtained as 
a union of a chain of subgroups, 
$$G=G_0(v)<G_1(v)<\ldots <G_n(v)\ldots ,$$
where $G_{i+1}(v)=G_i(v)\ast _{v_i=v_{i+1}^{i+1}}<v_{i+1}>;$ here $v_0=v.$
In other words, $G(C,\bf Q)$ can be obtained from $G$ as a union of a
 countable sequence
of elementary extensions of centralizers of the type 
\begin{equation}
\label{(1)}E(H,v,m)=H\ast _{v=w^m}<w>, \end{equation}
where the subgroup $<v>$ is maximal abelian in $H$. If $X$ is a fixed set
of generators of $H$, then we will consider the set $X \cup\{w \}$ as a 
canonical set of generators for $E(H,v,n)$. The length function on $E(H,v,n),$
introduced below,
is associated with this set of generators.

A cyclically minimal element $v$ of a group $G$ is called a 
{\em primitive} element if it is not a proper power. 

For an arbitrary group $G$ and natural number $n\geq 2$ 
choose a set of elements
${\cal V}_n = \{v_1\ldots v_t \}$ satisfying the following condition
${\em (S_n)}$:
\begin{enumerate}
\item[1)] ${\cal V}_n$ consists of primitive elements of length not more than 
$n$;
\item[2)] no two centralizers in the set of centralizers ${\cal C}_n=
\{C(v)| v\in {\cal V}_n\}$ are conjugate in $G;$
\item[3)] the set ${\cal V}_n$ is maximal with respect to 
properties 1) and 2); i.e.,
any element of length not more than $n$ is conjugate to a power of some
 $v\in {\cal V}_n$.
\end{enumerate}
By definition, the group $G({\cal C}_n)$ is the union of the finite chain of 
groups
$$G<E(G,v_1,n)=G_1< E(G_1,v_2,n)=G_2<\ldots <E(G_{t-1},v_t,n)=G_t=
G({\cal C}_n);$$
thus $G({\cal C}_n)$ is obtained from $G$ by consecutive extensions 
of centralizers from ${\cal C}_n:$
\begin{equation}
\label{(2)}G({\cal C}_n)=(\ldots 
(G\ast_{v_1=w_1^n}<w_1>)\ast_{v_2=w_2^n}<w_2>)\ast
\ldots )\ast_{v_t=w_t^n}<w_t>). \end{equation}
Notice that this definition does not depend on the order of the elements in
${\cal C}_n.$
\begin{lm} 
\label{(23)}Let $G({\cal C}_n)$ be as above. Then there exists a set 
${\cal V}_{n+1}$ in $G({\cal C}_n)$ that satisfies the condition ${\em 
(S_{n+1})}$
and contains  $\{w_1,\ldots ,w_t\}.$
\end{lm}
{\em Proof:} The elements $w_1, \ldots,w_t$ are simple in $G({\cal C}_n)$ 
because
their length is equal to 1. They are pair-wise nonconjugate in $G({\cal C}_n)$.
Indeed, from the description of conjugate elements
in free products with amalgamation (\cite{mks}), one can derive 
the following assertion:
 Let $g$ be a cyclically reduced element in $E(H,v,n)$.
Suppose $g$ is conjugate in $E(H,v,n)$ to some element in $H$.
Then $g\in H.$ 

Applying this fact to the chain (2) one can verify the pair-wise
nonconjugacy of the $w_i$'s. This proves the lemma.

We now intend to exhibit $G^{\bf Q}$ as a union of a chain of 
groups 
$$G=T_{0}<T_1<T_2<\ldots < \bigcup _{n=0}^{\infty}T_n,$$
defined by induction on $n$ as follows. Assume  that the 
groups $T_i, i<n,$
and  the sets ${\cal V}_i \subset T_{i-1}$  have
already been constructed and satisfy the  condition 
${\em (S_i)}$ in $T_{i-1}.$ If ${\cal V}_n=\{v_1,\ldots,v_t\},$ then

$$T_n=(\ldots (T_{n-1}\ast_{v_1=w_1^n}<w_1>)\ast_{v_2=w_2^n}<w_2>)\ast
\ldots )\ast_{v_t=w_t^n}<w_t>);$$
i.e. $T_n=T_{n-1}({\cal C}_n),$  where ${\cal C}_{n}$ is the corresponding 
set
of centralizers. By the previous lemma there exists a subset
${\cal V}_{n+1} \subset T_n$ which satisfies the condition {\em $(S_{n+1})$}
in $T_n$ and contains $w_1,\ldots ,w_t.$

\begin{df} We will call a sequence of groups $H_1,H_2,\ldots $
{\em effective} if there is an algorithm wich allows one, for any $i,$ to
construct a finite representation of the group $H_i.$\end{df}

\begin{th} 
\label{(4)}
Let $G$ be a torsion-free hyperbolic group
and $$G=T_0<T_1<T_2<\ldots $$ the sequence of groups described above;
then \begin{enumerate}
\item $T_n$ is hyperbolic for any $n,$
\item $\{T_n\}_{n\in\bf N}$ is an effective sequence of groups,
\item $ \bigcup _{n=0}^{\infty}T_n=G^{\bf Q}.$
\end{enumerate}
\end{th}
{\em Proof:} 1. By our construction, $T_{n+1}$ can be obtained from $T_n$
by a finite sequence of extensions of centralizers of the type 
$E(H,v,m)=H\ast _{v=w^m}<w>,$ where the subgroup $<v>$ is  maximal abelian
in $H$.  Suppose $H$ is torsion-free hyperbolic,
then $<v>$  is conjugate separated in $H$ and, by Corollary 2, $E(H,v,m)$
is also hyperbolic. $T_0=G$ and $G$ is hyperbolic; hence, by induction, 
all $T_n$ are hyperbolic.
Moreover, the proof of Theorem 2 shows that the constant
of hyperbolicity $\delta (E(H,v,m))$ can be effectively found if we know
a finite presentation of $H.$ This means that for any $n$ we can find 
effectively the constant of hyperbolicity of the group $T_n$, if we know
such a constant for $G$.

2. Recall that $T_{n+1}=T_n({\cal C}_{n+1}); $ so, having a finite 
presentation for $T_n$ and an effective procedure to construct the 
set of elements ${\cal V}_n,$ we can effectively construct a finite 
presentation for $T_{n+1}$ (see the presentation (2) above). 

\begin{lm} \label{24}
If $G$ is a torsion-free hyperbolic group 
then there is an algorithm to construct
the set
${\cal V}_n$
for any number $n.$ 
\end{lm}
{\em Proof:}
To prove the assertion of the lemma we need the solvability in any
hyperbolic group of the 
word and   
conjugacy problems  as well as the power problem; i.e. the problem 
to decide if an element is a proper power in $G$. The word and conjugacy
problems are solved in \cite{G}.

To solve  the power problem we will use the quasigeodesic property
of powers in hyperbolic groups; namely,
there is a constant $\lambda >0$  such that
$|W^n|\geq\lambda n|W|$ for any cyclically minimal
word $W$ (see \cite{Ol}, Lemma 1.12 and
\cite{Ol1}, Lemma 27). Examining the proofs of these lemmas one sees 
that $\lambda$ can be found effectively as a function of $\delta$
and $|A|,$ where $A$ is the distinguished system of generators of $G$. 
In light of the solvability of the word and conjugacy problems,
it is possible to decide if an element is cyclically minimal,
and to enumerate all cyclically minimal elements in accordance with
increasing 
length. Note that if 
a cyclically minimal element $V$ is equal to $W^n$, then $W$
is cyclically minimal and $n|W|\leq \lambda |V|.$
Therefore, to determine if a cyclically minimal element $V$ is a proper 
power, we
enumerate all cyclically minimal elements $W$ such that
$|W|\leq \lambda |V|$  and for each $n$ such that 
$n|W|\leq \lambda |V|,$ we verify the equality $W^n=V.$
It shows that we can effectively list all simple elements of 
length not greater than $n$. Now, to create a set ${\cal V}_n$ one  only  
needs  
to list all simple elements of length not greater than $n,$ and delete
those which are conjugate to previous ones. It can be done effectively
because of the decidability of the conjugacy problem.

3. As we mentioned above, every  element in ${\cal V}=\bigcup {\cal V}_n$
has arbitrary roots in the union $T= \bigcup _{n=0}^{\infty}T_n.$
Moreover, every centralizer $C_T(v), v \in {\cal V},$ is isomorphic 
to the additive group of ${\bf Q}$, so it admits an action of ${\bf Q}$ 
satisfying the module axioms.
By our construction, every centralizer in $ T$ is conjugate to 
the centralizer of an element $v \in {\cal V}.$ Hence  
we can, through conjugation, induce 
an action of ${\bf Q}$ on all centralizers in $T$; i.e. on the group $T$ .
This action is defined unambiguosly, because different centralizers in $T$ 
have trivial intersection (see \cite{MyasExpo2}).   
Hence $ T$ is a ${\bf Q}$-group. Using
the universal property of free products with amalgamation,
one can prove that this group satisfies the universal property of a 
$\bf Q$-completion of $G.$ See \cite{MyasExpo2} for details.
The theorem has been proved.

Let us now discuss algorithmic problems over groups of type $G^{\bf Q}$.
Formally, ${\bf Q}$-groups can be considered as groups with operators from 
${\bf Q}$. This means that the language of ${\bf Q}$-groups contains group
multiplication and countably many operations $f_{\alpha}, \alpha \in {\bf Q}$
(here, by definition, $f_\alpha(g)=g^{\alpha}$ for any $g \in G$). There are 
free objects (free ${\bf Q}$-groups) in the variety of ${\bf Q}$-groups, so, 
as usual, one can consider presentations of ${\bf Q}$-groups in the variety
of ${\bf Q}$-groups. 
\begin{lm} Let $<X|R>$ be a presentation of a group $G$. Then the group
$G^{\bf Q}$ has the same presentation $<X|R>$ in the variety of ${\bf 
Q}$-groups.
\end{lm}
   The proof follows from the fact that the ${\bf Q}$-group with presentation 
$<X|R>$ and $G^{\bf Q}$ have the same universal property in the category of 
${\bf Q}$-groups.

Therefore, ${\bf Q}$-completions of hyperbolic groups have finite presentations
as ${\bf Q}$-groups. Let $A=\{a_1,\ldots,a_m \}$ be a generating set for $G$. 
Then the elements in $G^{\bf Q}$ can be represented 
in the language of $\bf Q$-groups 
by words in the alphabet
$A.$  We are now able to formulate naturally
the word problem, the conjugacy problem, the equations problem and 
the isomorphism problem for finitely presented ${\bf Q}$-groups. Furthermore, 
when
we say that some algorithm is applied to some elements or groups, this means 
that
the algorithm is applied to the corresponding ${\bf Q}$-words and 
${\bf Q}$-presentations. 
   We would like to reduce the algorithmic problems for $G^{\bf Q}$ to the
hyperbolic groups $T_n$. To this end we need to construct 
 an algorithm to determine for 
any element $g\in G^{\bf Q}$ the number $n(g)$ such that
$g\in T_{n(g)}.$ We will do it at the same time as we construct normal 
forms for elements from $G^{\bf Q}.$ These normal forms derive from the 
construction
of $G^{\bf Q},$  starting from $G,$ as a countable iteration of extensions by 
free 
products with 
amalgamation.
   First of all, we introduce normal forms for the elements in the group
$$E(H,v,m)=H\ast _{v=w^m}<w>,$$
which is the elementary extension of a centralizer $C(v)=<v>$ in $H$ by the 
adjunction of an $m^{\mbox{th}}$ root to $v$.
     It will be convenient to denote $w$ by $v^{\frac{1}{m}}$. Let 
$S_m=\{\frac{k}{m}| k=1,\ldots, m-1\}$.  With this notation, the 
set $\{v^s| s\in S_m \}$ is a system of representatives in 
$<w>$ for the cosets of $<v>$. 
\begin{df}

A sequence of elements $(h_1, v^{s_1}, h_2, \ldots, v^{s_n}, h_{n+1}),$ 
where $h_i \in H, s_i \in S_m,$ is a {\em semicanonical form} of an
 element $g \in E(H,v,m)$ if 
\begin{equation}
\label{4}
g=h_1v^{s_1} \cdots v^{s_n}h_{n+1} 
\end{equation}
and  $h_i\not \in <v>, i=2, \ldots, n$.
\end{df}

\begin{lm}
Any two semicanonical forms of $g$ can be transformed from one to another
by a finite sequence of commutations of the form $v^sv^t=v^tv^s$,
where $t\in {\bf Z}, s\in S_m.$ 
\end{lm}
The proof follows from the definition of reduced forms for elements of 
a free product with 
amalgamation (see \cite{MyasExpo2} for details).

    Taking fixed right coset representatives of the subgroup $<v>$  in $H$ 
as the elements 
$h_i, i=2,\ldots,n$ in 
(\ref{4}) we obtain the notion of {\it canonical form } of $g$.    

Suppose now that the notions of canonical and semicanonical
forms of elements of the group $H$ have already 
been introduced. One can then extend them to the group $E(H,v,m)$.\par
We will say that the semicanonical (resp. canonical) forms of an element
$$
g=h_1v^{s_1} \cdots v^{s_n}h_{n+1}
$$
in $E(H,v,m)$ {\it agree with those on $H$}, iff the elements
$h_1, \ldots, h_{n+1}, v$ are in semicanonical (resp. canonical)
form in the group $H$. By definition an element $h \in H$ has the
same semicanonical (resp. canonical) form in $E(H,v,m)$ as in $H$.

Now let us consider  the union of a chain of groups
$$
G = G_0 < G_1 < \ldots < G_n< \ldots < \bigcup_{n=0}^{\infty} G_n, 
$$
where every $G_{n + 1}$ is obtained from $G_n$
by an elementary extension of a centralizer (i.e. it is an extension of the type
 $E(H,v_{n+1}, m_{n+1})$ and the notions of
semicanonical and canonical forms on $G_{n +1}$
are compatible with those on $G_n$. The forms on the terms of the chain
induce corresponding well-defined 
forms on the resulting group $\bigcup_{n=0}^{\infty} G_n .$
 Note that the described forms 
depend on the chosen elements  $ v_n, n\in{\bf N},$ and numbers  $m_n, n \in
{\bf N}$.

  By our construction, the ${\bf Q}$-completion $G^{\bf Q}$ of a torsion-free
hyperbolic group $G$ is the  union of the chain
$$
G=T_0 < T_1< \ldots < T_n< \ldots
,$$ 
where $T_{n+1}$ is obtained from $T_n$ by a finite sequence of elementary 
extensions of centralizers. Hence, to introduce the semicanonical 
and canonical forms on 
$G^{\bf Q}$ it is enough to introduce them on $G$. 
\begin{df} Let $A$ be the  generating set of $G.$ 
The  canonical (as well as the  semicanonical) 
form of an element $g\in G$ 
 is an $A$-word of
miminal length representing 
$g.$ 
The corresponding forms on $G^{\bf Q}$ (as described above) are 
called {\em induced semicanonical (canonical) forms}.
These forms depend on the chosen sets 
${\cal V}_n, n \in {\bf N},$
of elements $v_i.$
 
\end{df}
Let us suppose that some fixed sets of elements 
${\cal V}_n, n \in {\bf N},$ have been chosen.
\begin{lm} Every element $g \in G^{\bf Q}$ has a semicanonical form of 
the type $$(x_1, v_1^{s_1}, x_2,\ldots , v_m^{s_m}, x_m),$$ where
$x_i \in G, v_i \in {\bigcup {\cal V}_n}, s_i \in \bigcup S_n$.
\end{lm}
  The proof by induction is straightforward.

By definition 20, a semicanonical form of an element $g\in G^{\bf Q}$ 
is a sequence
of group elements $(h_1, v^{s_1}, h_2, \ldots, v^{s_n}, h_{n+1})$. 
When discussing 
algorithmic problems in the variety of ${\bf Q}$-groups, we will also consider 
semicanonical forms as sequences
of ${\bf Q}$-words representing the corresponding elements. Moreover, if   
the sets of words ${\cal V}_n, n \in {\bf N},$ are fixed, then the words  
representing elements $v_i$ in semicanonical form must be fixed 
words from ${\cal V}_n$, and not arbitrary words representing $v_i$ in 
$G^{\bf Q}$.

\begin{lm}
\label{(28)}There is an algorithm which  for every element $g \in G^{\bf Q},$ 
given in 
 semicanonical form, computes a number $n=n(g)$ such that $g \in T_n$.
\end{lm}
{\em Proof:} 
 According to lemma \ref{24}, there exists an algorithm which, for each $n,$ 
lists
the elements of the set ${\cal V}_n$ (i.e. it lists $\bf Q$-words in the  
alphabet $A$
representing these elements). For a given semicanonical form  
$(x_1, v_1^{s_1}, x_2,\ldots, v_m^{s_m}, x_m)$ one can effectively find for 
each $v_i$ a number $n_i$ such that $v_i \in {\cal V}_{n_i}$. By our 
construction, this semicanonical form represents an element from the subgroup
$T_n$, where $n$ is the maximal number in the set consisting of 
all $n_i$'s and the denominators of all $s_i$'s.

\begin{th}
\label{(5)}
Let $G$ be a torsion-free hyperbolic group, then there are algorithms that,
for any element $g \in G^{\bf Q},$ produce its semicanonical and canonical 
forms
with respect to some fixed family of sets ${\cal V}_n, n \in {\bf N}.$
\end{th}
{\em Proof:} Let us fix an arbitrary family of sets 
${\cal V}_n, n \in {\bf N}$ 
(which can be 
computed by some algorithm).
 For any ${\bf Q}$-word $W$ in the alphabet $A,$ representing some element 
$g \in G^{\bf Q},$ one needs to construct effectively the  
canonical form of $g$. We will argue by induction on the depth of the
word $W$. The {\em depth} of $W$ is a positive integer $d(W),$ 
defined by induction: $d(a)=1$ for any $a \in A$; 
$d(W_1W_2)=\max\{d(W_1),d(W_2)\}$ for any ${\bf Q}$-words $W_1, W_2$;
$d(W^{\alpha})=d(W)+1,$ where $\alpha \in {\bf Q}\setminus {\bf Z}$. It is easy 
to 
see that there exists an algorithm which 
for any {\bf Q}-word $W$ calculates $d(W)$.
 
Let $d(W)=1,$ then $W$ represents an element from the initial group $G$.
In light of the decidability of the word problem in $G,$ one 
can effectively construct
a canonical form of $g$, i.e. a minimal word in the 
alphabet $A$  representing
the element $g$.
   
Let $d(W) >1$. Then $W=W_1^{r_1}\ldots W_h^{r_h}$, where $d(W_i)< d(W), 
r_i \in  {\bf Q},$ and the ${\bf Q}$-words $W_1, \ldots, W_k$ can be found 
effectively from $W$. By induction, we can effectively determine 
for any $i$ the 
semicanonical
form of $W_i$ and therefore determine a number $n$ such that $W_i \in 
T_n$. By the definition of the sets ${\cal V}_n, n \in {\bf N},$
if an integer $k$ is greater then the length of all $W_i$ in the generators of 
the
group $T_n,$ then  any $W_i$
is conjugate to some power of an element $v_i$ from the set ${\cal V}_{n+k}$.
By lemma \ref{24} the elements from the set ${\cal V}_{n+k}$ can be listed 
effectively.
So, looking through all words $x$ (for example, according to increasing 
lengths) one can find effectively the words $v_i \in {\cal V}_{n+k}$, 
some integers $l_i$
and words $x_i$ such that $W_i=x_i^{-1}V_i^{l_i}x_i$ in the group $T_{n+k}$.
Using axiom 3) from the definition of an $A$-group, the word $W$ 
represents 
in $T_{n+k}$ the same element as the word 

$$ x_1^{-1}V_1^{l_1r_1}x_1 \ldots x_k^{-1}V_k^{l_kr_k}x_k,$$
and the latter word is constructed effectively from $W$. So it is left only
to transform it into semicanonical form. The procedure for this is the 
following:
find equal neighbours $V_i=V_{i+1}=V$ (if any) and verify if the word
$x_ix_{i+1}^{-1} $ between them belongs to the cyclic subgroup $<V>$ 
(it can be done effectively, as explaned before). If $x_ix_{i+1}^{-1}= V^t $
for some integer $t,$ then  replace the subword
$V^{r_il_i}x_ix_{i+1}^{-1}V^{r_{i+1}l_{i+1}}$     
by the word 
$V^q$, where $q= r_il_i+r_{i+1}l_{i+1}+t$. And if $q=m+{\frac{s_1}{s_2}}$, 
where $m$ is an integer and $0<s_1<s_2$, then replace $V^q$ by the word
$V^mV^{\frac{s_1}{s_2}}$. The resulting word represents the same element as 
the original word, but the number of $V_i$'s in it is less than in
the original one. Arguing by induction 
we complete the process. It is not hard to see that the resulting word will
be in semicanonical form.
The theorem 
has been proved.

\begin{th} 
\label{wcp}
Let $G$ be a torsion-free hyperbolic group, then the
word and conjugacy problems are solvable in $G^{\bf Q}.$
\end{th}
{\it Proof:} Let $W$ and $V$ be arbitrary {\bf Q}-words in the alphabet $A$.
By Theorem \ref{(5)} and Lemma \ref{(28)}, one 
can effectively find a number $n=n(W,V)$
such that the words $W, V$ represent some elements $ g,h \in T_n$ of  
$G^{\bf Q}.$ By Theorem \ref{(4)} the group $T_n$ is hyperbolic and 
there is an algorithm to construct  a finite 
presentation of $T_n$. As was noticed in Lemma~\ref{(23)},  if 
$g,h$ are conjugate
in $G^{\bf Q}$ then they are conjugate in $T_n$. So we have reduced the 
word and conjugacy problems to the hyperbolic case. But in hyperbolic groups 
the above-mentioned problems are decidable. 
This completes the proof of the Theorem.


\begin{thebibliography}{10}

\bibitem{baum1}
G Baumslag.
\newblock On free ${\bf Q}$-group.
\newblock {\it Commun. on Pure and Appl. Mathematics}, 18:25--30, 1965.

\bibitem{bau}
G. Baumslag.
\newblock Some aspects of groups with unique roots.
\newblock {\it Acta Math.}, 104:217--303, 1960.

\bibitem{BGSS1}
G. Baumslag, S.M. Gersten, M. Shapiro, and H. Short.
\newblock Automatic groups and amalgams.
\newblock 1990.
\newblock Preprint, Ohio State University.

\bibitem{BGSS}
G. Baumslag, S.M. Gersten, M. Shapiro, and H. Short.
\newblock Automatic groups and amalgams - a survey.
\newblock In {\it Algorithms and classification in combinatorial group theory,
  MSRI, 23}, pages~179--194, 1991.

\bibitem{BF}
M. Bestvina and M Feighn.
\newblock A combination theorem for negatively curved groups.
\newblock {\it J. Diff. Geom.}, 35:85--101, 1992.

\bibitem{GS}
S.M. Gersten and H.B. Short.
\newblock Rational subgroups of biatomatic groups.
\newblock {\it Annals of Mathematics}, 134:125--158, 1991.

\bibitem{Rita}
R. Gitic.
\newblock On combination theorem for negatively curved groups.
\newblock 1995.
\newblock to appear in IJAC.

\bibitem{G}
M. Gromov.
\newblock {\it Hyperbolic groups}, pages~75--263.
\newblock Springer, Berlin, 1987.

\bibitem{IK}
I. Kapovich.
\newblock On a theorem of {B.} {Baumslag}.
\newblock 1994.
\newblock Proceedings of the AMS Meeting, Brooklyn, New-York, to appear.

\bibitem{mks}
W. Magnus, A. Karras, and D. Solitar.
\newblock {\it Combinatorial group theory}.
\newblock Interscience publ., 1966.

\bibitem{OM}
K.V. Mikhajlovskii and A.Yu. Ol'shanskii.
\newblock Some constructions relating to hyperbolic groups.
\newblock 1994.

\bibitem{MyasExpo2}
A.~G. Myasnikov and V.~N. Remeslennikov.
\newblock Exponential groups 2: extension of centralizers and tensor completion
  of csa-groups.
\newblock 1994.
\newblock Accepted to the IJAC.

\bibitem{Ol1}
A.Yu. Olshanski.
\newblock Periodic factor groups of hyperbolic groups.
\newblock {\it Math. USSR Sbornik}, 72(2):512--541, 1992.

\bibitem{Ol}
A.~Yu. Ol'shanskii.
\newblock On residualing homomorphisms and {G-subgroups} of hyperbolic groups.
\newblock {\it Int. J. Algebra and Computation}, 3(4):365--409, 1993.

\bibitem{PP}
P. Papasoglu.
\newblock {\it Geometric methods in group theory}.
\newblock PhD thesis, Columbia University, 1993.

\bibitem{Sh}
M.D. Shapiro.
\newblock Automatic structure and graphs of groups.
\newblock In {\it Topology'90 (Columbus, OH,90)}, pages~335--380, Math. Res.
  Inst. Publ., 1 de Gruyer, Berlin, 1992.

\end{thebibliography}
\bibliographystyle{/pctex/texbib/plain}

\end{document}